\documentclass[reqno]{amsart}
 \usepackage{latexsym}
 \usepackage{amsbsy}
 \usepackage{amssymb}
 \usepackage{amsfonts}
 \usepackage{amsmath}
 \usepackage[latin1]{inputenc}
 \usepackage{graphicx}
 \usepackage{epsfig}
 \input{amssym.def}
 \input{amssym}
 \newtheorem{thm}{Theorem}[section]
 \newtheorem{defin}[thm]{Definition}
 \newtheorem{lem}[thm]{Lemma}
 \newtheorem{prop}[thm]{Proposition}
 \newtheorem{cor}[thm]{Corollary}
 \newtheorem{rem}[thm]{Remark}
 \newtheorem{ex}[thm]{Example}

 \newcommand{\bthm}{\begin{thm}}
 \newcommand{\ethm}{\end{thm}}
 \newcommand{\bd}{\begin{defin}}
 \newcommand{\ed}{\end{defin}}
 \newcommand{\blem}{\begin{lem}}
 \newcommand{\elem}{\end{lem}}
 \newcommand{\bcor}{\begin{cor}}
 \newcommand{\ecor}{\end{cor}}
 \newcommand{\bprop}{\begin{prop}}
 \newcommand{\eprop}{\end{prop}}
 \newcommand{\brem}{\begin{rem} \rm}
 \newcommand{\erem}{\end{rem}}
 \newcommand{\bex}{\begin{ex} \rm}
 \newcommand{\eex}{\end{ex}}

 \newcommand{\pr}{\noindent{\bf Proof. }}
 \newcommand{\ep}{\nolinebreak{\hspace*{\fill}$\Box$ \vspace*{0.25cm}}}

 \newcommand{\ds}{\displaystyle}

 \def\dj{d\kern-0.4em\char"16\kern-0.1em}
 \def\Dj{\mbox{\raise0.3ex\hbox{-}\kern-0.4em D}}

 \newcommand{\beq}{ \begin{equation} }
 \newcommand{\eeq}{\end{equation} }
 \newcommand{\bea}{\begin{eqnarray}}
 \newcommand{\eea}{\end{eqnarray}}
 \newcommand{\beast}{\begin{eqnarray*}}
 \newcommand{\eeast}{\end{eqnarray*}}
 \newcommand{\ba}{\begin{array}}
 \newcommand{\ea}{\end{array}}
 \def\bi{\begin{itemize}}
 \def\ei{\end{itemize}}
 \newcommand{\bc}{\begin{center}}
 \newcommand{\ec}{\end{center}}

 \newcommand{\eps}{\varepsilon}
 \newcommand{\mR}{\mathbb R}

\begin{document}
 \title{Delta shock wave interactions via wave front tracking method}
 \author{Neboj\v sa Dedovi\'c and Marko Nedeljkov}
 \maketitle

 \begin{abstract}
 In this paper we discuss delta shock interaction problem for a
 pressureless gas dynamics system with two different ways of approaching the
 subject. The
 first one is by using shadow wave solution concept.
 The result of two delta shock interactions is
 delta shock with non-constant speed in a general case. The second
 one is by perturbing the system with a small pressure term.
 The obtained perturbed system is strictly hyperbolic and its Riemann problem is solvable.
 We compare a limit of a numerical wave front tracking results as
 small pressure term vanishes with the shadow wave solution.

 \end{abstract}

 \section{Introduction}
 \thispagestyle{empty}

 Consider the one-dimensional Euler gas dynamics system given by
 \beq
 \label{modsys}
 \begin{array}{rcl}
 \partial_t\rho + \partial_x(\rho u) & = & 0\\
 \ds \partial_t(\rho u)+\partial_x(\rho u^2+p(\eps,\rho)) & = & 0,
 \end{array}
 \eeq

 where $\rho$ is the density, $m=\rho\,u$ is the momentum,
 $p(\eps,\rho)=\eps\,p_0(\rho)$ is the scalar pressure, $\eps<<1$
 and $p_0(\rho)=\rho^\gamma/\gamma\,$. Taking $\eps\to 0$ in (\ref{modsys}), we obtain the
 pressureless gas dynamics model (PGD model in the
 rest of the paper), also called sticky particles model (in \cite{ERS})
 \beq
 \label{PGD}
 \begin{array}{rcl}
 \partial_t\rho + \partial_x(\rho u) & = & 0\\
 \ds \partial_t(\rho u)+\partial_x(\rho u^2) & = & 0,\; (x,t)\in
 \mR\times\mR_+\,.
 \end{array}
 \eeq

 System (\ref{modsys}) can be considered as a perturbation of
 system (\ref{PGD}) which is weakly hyperbolic with a double eigenvalue $\lambda_1=\lambda_2=u$.
 All entropy pairs $(\eta,q)$ with a semiconvex function $\eta$
 are given by $\eta:=\rho S(u)$, $q:=\rho\,u S(u)$, where $S''\geq 0$
 (the entropy function $\eta$ is semi-convex with respect to the
 variable $(\rho,\rho u)$).
 The Riemann problem
 \beq\label{idPGD}
 \rho(x,0)=\left\{
 \begin{array}{cc}
 \rho_0 , & x<0\,, \\
 \rho_1 , & x>0\,,
 \end{array}
 \right.\;,\hspace*{1cm}
 u(x,0)=\left\{
 \begin{array}{cc}
 u_0 , & x<0\,, \\
 u_1 , & x>0\,,
 \end{array}
 \right.
 \eeq

 has a classical entropy solution consisting of two contact
 discontinuities connected with the vacuum state ($\rho=0$) if
 $u_{0}\leq u_{1}$:
 \begin{equation*}
 (\rho(x,t),u(x,t))=\begin{cases}
 (\rho_{0},u_{0}), & x<u_{0}t, \\
 (0,\psi(x/t)), & u_{0}t<x<u_{1}t, \\
 (\rho_{1},u_{1}), & x>u_{1}t,
 \end{cases}
 \end{equation*}

 where $\psi(y)=y$. We are now turning to the case $u_{0}>u_{1}$ when there is no
 classical solution to the Riemann problem (\ref{PGD}, \ref{idPGD}).

 Throughout this paper, the following constants will be fixed:
 \beq\label{definP}
 \gamma=1+2\eps,\;0<\eps<\frac12,\;
 \kappa=\frac{\sqrt\eps}{\sqrt\gamma}\;\;{\rm and }\;\;p=\kappa^2\rho^\gamma\;.
 \eeq

 \section{Elementary waves of the perturbed system}\label{rwswcurves}

 The eigenvalues of system (\ref{modsys}) are
 \beq
 \label{charvalues}
 \ba{l}
 \ds \lambda_1 =u-\kappa \sqrt{\gamma} \rho ^{\frac{\gamma-1}{2}},\\
 \\
 \ds \lambda_2 =u + \kappa \sqrt{\gamma} \rho
 ^{\frac{\gamma-1}{2}},
 \ea
 \eeq

 and the corresponding eigenvectors are
 \beq
 \label{charvectors}
 \ba{l}
 \ds r_1=(-1,-u+\kappa \sqrt{\gamma} \rho
 ^{\frac{\gamma-1}{2}})^T,\\
 \\
 \ds r_2=(1,u+\kappa \sqrt{\gamma} \rho
 ^{\frac{\gamma-1}{2}})^T\,.
 \ea
 \eeq

 We have chosen an orientation such that $\nabla\lambda_i\cdot
 r_i>0,\;i=1,2$,
 since both fields are genuinely nonlinear. The corresponding Riemann invariants
 of system (\ref{modsys}) are
 \beq
 \label{riemann inv}
 \ba{c}
 s=u+\frac{\kappa\sqrt{\gamma}}{\eps}(\rho^{\eps}-1):\mbox{
 1-invariant, and}\\
 r=u-\frac{\kappa\sqrt{\gamma}}{\eps}(\rho^{\eps}-1):\mbox{ 2-invariant}\,.
 \ea
 \eeq

 The rarefaction curves through the point $(\rho_0, u_0)$ are given by
 \beq
 \label{rarefaction waves}
 \ba{c}
 u-u_0=-\frac{\kappa\sqrt{\gamma}}{\eps}(\rho^{\eps}-\rho^{\eps}_0)\;,\;\;0\leq\rho\leq\rho_0:
 \mbox{ 1-rarefaction curve},\\
 u-u_0=\;\;\;\frac{\kappa\sqrt{\gamma}}{\eps}(\rho^{\eps}-\rho^{\eps}_0)\;,\;\;\rho
 \geq\rho_0: \mbox{ 2-rarefaction curve},
 \ea
 \eeq

 while the shock curves through the point $(\rho_0, u_0)$ are given by
 \beq
 \label{shock wave eq1 new}
 u-u_0=-\kappa\sqrt{\frac{\rho^{\gamma}-\rho^{\gamma}_0}{\rho_0\rho(\rho-\rho_0)}}\;(\rho-\rho_0),\;\;\;\rho>\rho_0:
 \mbox{ 1-shock curve},
 \eeq

 and
 \beq
 \label{shock wave eq2 new}
 u-u_0=\kappa\sqrt{\frac{\rho^{\gamma}-\rho^{\gamma}_0}{\rho_0\rho(\rho-\rho_0)}}\;(\rho-\rho_0),\;\;\;0<\rho<\rho_0:
 \mbox{ 2-shock curve}.
 \eeq

 With the Riemann invariants, shock curves starting from the point $(r_0,s_0)$
 are
 \beq
 \label{s1 new}
 S_1:\;\;\; \left\{
 \ba{l}
 \ds r_0-r=\kappa\rho_0^{\eps}\left(\sqrt{\frac{(\alpha-1)(\alpha^\gamma-1)}{\alpha}}+\sqrt\gamma\frac{\alpha^{\eps}-1}{\eps}\right),\\
 \ds s_0-s=\kappa\rho_0^{\eps}\left(\sqrt{\frac{(\alpha-1)(\alpha^\gamma-1)}{\alpha}}-\sqrt\gamma\frac{\alpha^{\eps}-1}{\eps}\right),
 \ea\right.
 \eeq

 where $r_0=r(\rho_0,u_0)$, $s_0=s(\rho_0,u_0)$ and $\alpha=\rho/\rho_0\geq 1$, and
 \beq
 \label{s2 new}
 S_2:\;\;\; \left\{
 \ba{l}
 \ds s_0-s=\kappa\rho_0^{\eps}\left(\sqrt{\frac{(1-\alpha)(1-\alpha^\gamma)}{\alpha}}+\sqrt\gamma\frac{1-\alpha^{\eps}}{\eps}\right),\\
 \ds r_0-r=\kappa\rho_0^{\eps}\left(\sqrt{\frac{(1-\alpha)(1-\alpha^\gamma)}{\alpha}}-\sqrt\gamma\frac{1-\alpha^{\eps}}{\eps}\right),
 \ea\right.
 \eeq

 where $r_0=r(\rho_0,u_0)$, $s_0=s(\rho_0,u_0)$ and
 $0<\alpha=\rho/\rho_0\leq 1$. The corresponding rarefaction curves are given by
 \beq
 \label{r1 new}
 R_1:\;\;\; r\geq r_0,\; s=s_0,
 \eeq

 and
 \beq
 \label{r2 new}
 R_2:\;\;\; s\geq s_0,\; r=r_0.
 \eeq

 It is clear that from (\ref{s1 new}, \ref{s2 new})
 we have that $r_0-r\geq s_0-s$ holds for $S_1$ curve and $s_0-s\geq r_0-r$ holds for $S_2$ curve,
 respectively.\\

 The Riemann problem for system (\ref{modsys}) with initial data
 (\ref{idPGD}) was solved by Riemann \cite{BR}, and the result is
 summarized in the following theorem (the proof can be found in
 Courant-Friedrichs \cite{RCKF} and Smoller \cite{JS}).

 \bthm\label{essentthm}\cite{FA}
 Consider system (\ref{modsys}) with initial data (\ref{idPGD}).
 Suppose that $u_1-u_0<\frac{\kappa\sqrt{\gamma}}{\eps}(\rho^{\eps}_1+\rho^{\eps}_0)$,
 or equivalently $s_0-r_1>-\frac{2\kappa\sqrt{\gamma}}{\eps}$. Then
 there exists a unique solution composed of constant states $(\rho_0,u_0)=(r_0,s_0)$,
 $(\rho_m,u_m)=(r_m,s_m)$ and $(\rho_1,u_1)=(r_1,s_1)$ separated by centered rarefaction or shock waves
 satisfying the following estimates:
 \beq
 \label{wz app}
 \ba{l}
 r(x,t)=r(\rho(x,t),u(x,t))\geq \min\{r_0,r_1\},\\
 s(x,t)=s(\rho(x,t),u(x,t))\leq \max\{s_0,s_1\}.
 \ea
 \eeq

 The {\it amplitude} of the waves is denoted by
 \beq
 \label{amplitude}
 \ba{c}
 \beta:=r_m-r_0\;\;: \mbox{ amplitude of an 1-wave},\\
 \chi:=s_1-s_m\;\;: \mbox{ amplitude of a 2-wave}.
 \ea
 \eeq

 Here $\beta, \chi\geq 0$ for centered rarefaction waves and $\beta,
 \chi<0$ for shock waves; absolute values $|\beta|$, $|\chi|$
 are called {\it strengths} of $\beta$ and $\chi$, respectively.
 \ethm

 We shall use that notation throughout the rest of the paper.

 \section{Local Interactions Estimates}

 Our first task is to obtain a sharp estimate of wave strengths with
 respect to $\eps$ as much as possible. In order to do that, we
 shall present some assertions from \cite{TNJS} together with
 modified proofs, since certain changes in estimates will be useful for our investigation.

 \bthm\label{helpthm}\cite{TNJS} The shock curve $S_1$ starting at the point $(r_0,s_0)$
 is given by
 \beq\label{g1fun}
 \ds
 s_0-s=g_1(r_0-r,\rho_0)=\int_0^{r_0-r}\;h_1(\alpha)|_{\alpha=\alpha_1(\beta/\kappa\rho^{\eps}_0)}\;d\beta,\;\;r<r_0,
 \eeq

 where $0\leq g'_1(\beta,\rho_0)<1$ and $g''_1(\beta,\rho_0)\geq 0$\footnote{The primes denote differentiation with respect to the
 first argument.}. The shock curve $S_2$ starting at the point
 $(r_0,s_0)$ is
 \beq
 \ds
 r_0-r=g_2(s_0-s,\rho_0)=\int_0^{s_0-s}\;h_2(\alpha)|_{\alpha=\alpha_2(\chi/\kappa\rho^{\eps}_0)}\;d\chi,\;\;s<s_0,
 \eeq

 where $0\leq g'_2(\chi,\rho_0)<1$ and $g''_2(\chi,\rho_0)\geq 0$.
 \ethm

 \pr We shall repeat the proof from \cite{TNJS} in order to fix the
 notation for the rest of the paper. Relation $s_0-s=g_1(r_0-r,\rho_0)$
 implies
 \beq
 \label{h1a}
 \frac{\partial (s_0-s)}{\partial \alpha}=\frac{\partial g_1(r_0-r,\rho_0)}{\partial
 (r_0-r)}\cdot \frac{\partial (r_0-r)}{\partial \alpha},\;{\rm so
 }\;\;
 \frac{\partial (s_0-s)/\partial \alpha}{\partial (r_0-r)/\partial
 \alpha}=g'_1(\beta,\rho_0)\,.
 \eeq

 If
 $$
 h_1(\alpha)=\frac{\partial (s_0-s)/\partial \alpha}{\partial (r_0-r)/\partial
 \alpha}\,,
 $$

 then one can easily see that
 \beq
 \label{h1}
 \ds h_1(\alpha)=\left(\frac{Y-1}{Y+1}\right)^2\;\;\;\mbox{ with }\;\;
 Y=\sqrt{\frac{\gamma\alpha^\gamma(\alpha-1)}{\alpha^\gamma-1}}\;,\;\;\mbox{ for }\;\;\alpha>1\;.
 \eeq

 From the first equation in (\ref{s1 new}), we have
 \beq
 \label{f}
 \frac{\beta}{\kappa\rho^{\eps}_0}=\sqrt{\frac{(\alpha-1)(\alpha^\gamma-1)}{\alpha}}+\sqrt\gamma\;\frac{\alpha^{\eps}-1}{\eps}=:
 f(\alpha)\;.
 \eeq

 Therefore
 $$
 \ba{rl}
 \ds f'(\alpha)> &\ds \frac{1}{2}\sqrt\frac{\alpha}{(\alpha-1)(\alpha^\gamma-1)}
 \cdot\frac{\alpha^\gamma-1}{\alpha^2}+\sqrt\gamma\;\alpha^{\eps-1}>0
 \ea
 $$

 since $\alpha^\gamma>1$ for $\alpha>1$ and $\gamma>1$.
 Using the fact that $f'(\alpha)>0$ and (\ref{f}) the Implicit Function Theorem yields that there exists
 $\alpha=\alpha_1(\beta/\kappa\rho^{\eps}_0)$ such that
 \beq
 \label{h1 new}
 g_1(r_0-r,\rho_0)=\int_0^{r_0-r}\;h_1(\alpha)|_{\alpha=\alpha_1(\beta/\kappa\rho_0^\eps)}\;d\beta\;.
 \eeq

 Since $g'_1(\beta,\rho_0)=h_1(\alpha)$,
 $\ds g''_1(\beta,\rho_0)=h'_1(\alpha)\cdot\frac{d\alpha}{d\beta}$ and
 $\ds\frac{d\beta}{d\alpha}=\kappa\rho_0^\eps\cdot f'(\alpha)>0$
 it remains to prove that $0\leq h_1(\alpha)<1$ and $0\leq
 h'_1(\alpha)$. From
 (\ref{h1}) we have
 $$
 0\leq
 h_1(\alpha)=\left(\frac{Y-1}{Y+1}\right)^2<\left(\frac{Y+1}{Y+1}\right)^2=1\;,
 $$

 and
 \beq\label{h1prim}
 0\leq h'_1(\alpha)=4\cdot\frac{Y-1}{(Y+1)^3}\cdot Y'\;,
 \eeq

 since $Y\geq 1$ and $Y'\geq 0$. The
 second part of the
 theorem can be proved using the same technique.
 \ep

 \blem\label{lemanovo} Let $\rho_0<\rho_1$ and $\beta/\kappa\rho^{\eps}_1<\theta<\beta/\kappa\rho^{\eps}_0$.
 Then
 \beq
 \frac{d\alpha}{d\theta}=\frac{1}{f'(\alpha)}=\frac{2Y}{\sqrt{\gamma}\alpha^\frac{\gamma-3}{2}(1+Y)^2}\,.
 \eeq
 \elem

 We would need an estimate of the difference of Riemann invariants across
 two shock waves which is more precise than the one in
 \cite{TNJS}. It is provided by the following theorem.

 \bthm\label{diffrieminv} Let $0< \eps < \frac12$, $s_0<s_1$, and take two $S_1$ curves
 originating at the points $(r_0,s_0)=(\rho_0,u_0)$ and
 $(r_0,s_1)=(\rho_1,u_1)$, which are continued to the points $(r,s)$ and
 $(r,s_2)$, respectively. Then we have
 \beq
 \label{ocena invarijanti}
 0\leq (s_0-s)-(s_1-s_2)\leq C_*\,\sqrt\eps\,(r_0-r)\,(s_1-s_0)\,,
 \eeq

 where $C_*$ is a constant independent of $\eps$, $\rho_0$ and $\rho_1$\,.
 \ethm

\begin{figure}[htbp]
 \centerline{
 \epsfig{figure=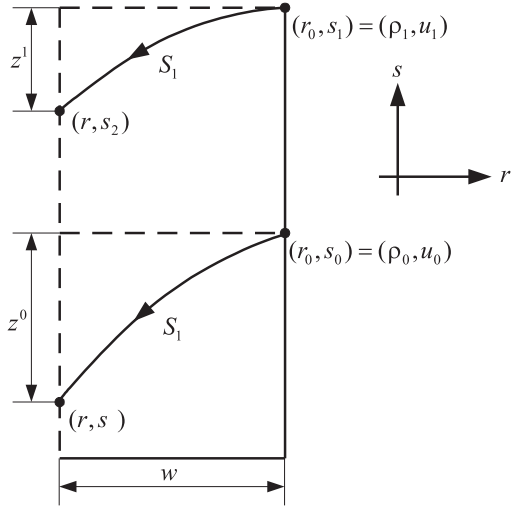,width=7cm
 }} \epsfxsize=1\textwidth
 \caption{Two 1-shock wave curves in $r-s$ plane.\label{s1s1}}
 \end{figure}

 \pr Let $z^0=s_0-s$, $z^1=s_1-s_2$ and $w=r_0-r$ (look at the
 diagram shown in Figure \ref{s1s1}).

 %
 %
 %
 %

 By Theorem \ref{helpthm} and the
 Mean Value Theorem we know that for $\rho_1>\rho_0$, there exists $\theta$
 such that
 \beq
 \label{z1z0}
 \ba{rl}
 z^0-z^1=&\ds\int_0^w
 \frac{dh_1(\alpha)}{d\alpha}\Big|_{\alpha=\alpha(\theta)}\cdot\alpha'(\theta)
 \left(\frac{\beta}{\kappa\rho_0^{\eps}}-\frac{\beta}{\kappa\rho_1^{\eps}}\right)\;d\beta\;,
 \ea
 \eeq

 where $\theta\in\left(\frac{\beta}{\kappa\rho^{\eps}_1},\frac{\beta}{\kappa\rho^{\eps}_0}\right)$.
 The definitions of $h_1$ and
 $\alpha$ imply
 $$
 \frac{dh_1(\alpha)}{d\alpha}\geq 0,\;\;\;
 \frac{d\alpha(\theta)}{d\theta}\geq 0\;\;\;\mbox{ and }\;\;\;
 \frac{\beta}{\kappa\rho_0^{\eps}}-\frac{\beta}{\kappa\rho_1^{\eps}}\geq
 0
 $$

 so $z^0-z^1\geq 0$ for $\rho_1>\rho_0$. We need to estimate the
 integrand in (\ref{z1z0}). By (\ref{h1}), we have
 $$
 \ba{rl}
 \ds
 \frac{dh_1(\alpha)}{d\alpha}\cdot\frac{d\alpha(\theta)}{d\theta}
 \leq&
 \ds
 \frac{4(Y-1)(\gamma+1)\alpha^\frac{1-\gamma}{2}}{\sqrt\gamma\;(Y+1)^3}\;.
 \ea
 $$

 Thus,
 \beq\label{z1z0new}
 \ba{rl}
 z^0-z^1\leq &\ds \frac{4(\gamma+1)}{\kappa\sqrt\gamma\;\rho_0^\eps \rho_1^\eps}(\rho_1^{\eps}-\rho_0^{\eps})
 \int_0^w
 \beta\;\alpha^\frac{1-\gamma}{2}\Big|_{\alpha=\alpha(\theta)}\frac{Y-1}{(Y+1)^3}\;d\beta
 \ea
 \eeq

 and
 \beq
 \label{z1z0novo1}
 \ba{rl}
 z^0-z^1\leq&\ds \frac{4(\gamma+1)}{\kappa\sqrt\gamma\;\rho_0^\eps \rho_1^\eps}(\rho_1^{\eps}-\rho_0^{\eps})
 \int_0^w
 \beta\;\alpha^{-\frac{1+\gamma}{2}}\Big|_{\alpha=\alpha(\theta)}\;d\beta\;.
 \ea
 \eeq

 From Lemma \ref{lemanovo} we know that $d\alpha/d\theta>0$ for
 $\beta/\kappa\rho^{\eps}_1<\theta<\beta/\kappa\rho^{\eps}_0$.
 Hence,
 \beq
 \label{alfapomoc}
 \ds
 \alpha\left(\frac{\beta}{\kappa\rho_1^{\eps}}\right)\leq\alpha\left(\theta\right)\leq
 \alpha\left(\frac{\beta}{\kappa\rho_0^{\eps}}\right).
 \eeq

 Moreover,
 $$
 \ba{rl}
 \ds \frac{\beta}{\kappa\rho^{\eps}_1}=f(\alpha)
 \leq &\ds
 2\sqrt{\frac{(\alpha-1)\;\alpha^\gamma}{\alpha-1}}=2\alpha^{\gamma/2}\;,\mbox{
 for }
 \alpha=\alpha\left(\frac{\beta}{\kappa\rho_1^{\eps}}\right)\;.
 \ea
 $$

 At this point, we use a majorization of $z^0-z^1$ different from the one in \cite{TNJS} in order to obtain bounds
 for $C_*$ independent of $\eps$. By (\ref{alfapomoc}) and the above inequality we obtain
 \beq \label{betakon}
 \ba{l}
 \ds \frac{\beta}{2\kappa\rho^{\eps}_1}\leq \left(\alpha\left(\frac{\beta}{\kappa\rho_1^{\eps}}\right)\right)^{\gamma/2}
 \Rightarrow \left(\frac{\beta}{2\kappa\rho^{\eps}_1}\right)^{-\frac{\gamma+1}{\gamma}} \geq
 \left(\alpha\left(\frac{\beta}{\kappa\rho^{\eps}_1}\right)\right)^{-\frac{\gamma+1}{2}}\,,\;\;\rm{ so}\\
 \\
 \ds \min\left\{1,\left(\frac{\beta}{2\kappa\rho^{\eps}_1}\right)^{-\frac{\gamma+1}{\gamma}}\right\} \geq
 \left(\alpha\left(\frac{\beta}{\kappa\rho^{\eps}_1}\right)\right)^{-\frac{\gamma+1}{2}}\geq
 \left(\alpha(\theta)\right)^{-\frac{\gamma+1}{2}}\;,
 \ea
 \eeq

 for $\alpha\left(\frac{\beta}{\kappa\rho_1^{\eps}}\right)\geq 1$ and $\gamma>1$.
 Since $d\alpha/d\theta>0$, it follows by (\ref{z1z0novo1}) that
 \beq
 \label{z1z0final}
 \ba{rl}
 z^0-z^1\leq&\ds \frac{4(\gamma+1)}{\kappa\sqrt\gamma\;\rho_0^\eps \rho_1^\eps}(\rho_1^{\eps}-\rho_0^{\eps})
 \int_0^w
 \beta\cdot\min\left\{1,\left(\frac{\beta}{2\kappa\rho^{\eps}_1}\right)^{-\frac{\gamma+1}{\gamma}}\right\}\;d\beta\\
 &\\
 \leq&\ds \frac{8(\gamma+1)}{\sqrt\gamma\;\rho_0^\eps}\;(\rho_1^{\eps}-\rho_0^{\eps})
 \int_0^w
 \min\left\{\frac{\beta}{2\kappa\rho^{\eps}_1},
 \left(\frac{\beta}{2\kappa\rho^{\eps}_1}\right)^{-\frac{1}{\gamma}}\right\}\;d\beta\\
 &\\
 \leq&\ds \frac{8(\gamma+1)}{\sqrt\gamma\;\rho_0^\eps}\;(\rho_1^{\eps}-\rho_0^{\eps})\;w\,.
 \ea
 \eeq

 Using
 $$
 \rho^{\eps}_1-\rho^{\eps}_0=\frac{\eps}{2\kappa\sqrt\gamma}(s_1-s_0)
 $$

 and $\kappa\sqrt\gamma=\sqrt\eps$ together with (\ref{z1z0final}) we finally get
 \beq
 \label{z1z0finalnew11}
 \ba{rl}
 z^0-z^1\leq&\ds \frac{4(\gamma+1)}{\sqrt\gamma\;\rho_0^\eps}\cdot\frac{\eps}{\kappa\sqrt\gamma}\;(s_1-s_0)\;w\\
 &\\
 =&\ds
 \frac{4(\gamma+1)}{\sqrt\gamma\;\rho_0^\eps}\cdot\sqrt\eps\;(s_1-s_0)\;w\,.
 \ea
 \eeq

 Suppose that $\rho$ is the first component of the
 solution of the Riemann problem (\ref{modsys}, \ref{idPGD}) for $u_0>u_1$.
 Then Lemma 3.1 from \cite{GQCHL} yields that for small $\eps>0$, there exists
 $C>0$ independent of $\eps$,
 such that $\ds \rho\leq C/\eps$. Now, using (\ref{z1z0finalnew11}), if $\rho_0\sim 1/\eps$
 then $\rho_0^\eps\sim 1$ as $\eps\to 0$. For $\eps$ small enough we may write
 $\rho_0^\eps \geq C_1$. Thus, there exists a constant $C_*$ independent of $\kappa$, $\rho_0$, $\rho_1$ and $\eps$ such
 that
 \beq
 \label{z1z0finalrealy}
 z^0-z^1\leq C_*\,\sqrt\eps\;(s_1-s_0)\;w,
 \eeq

 holds. This completes the proof of the theorem.
 \ep

 The theorem that follows can be proved in the same way.
 \bthm\label{diffrieminv2} Let $0< \eps < \frac12$, $r_0>r_1$, and take two $S_2$ curves
 originating at the points $(r_0,s_0)=(\rho_0,u_0)$ and
 $(r_1,s_0)=(\rho_1,u_1)$, which are continued to the points $(r,s)$ and
 $(r_2,s)$, respectively. Then we have
 \beq
 \label{ocena invarijanti2}
 0\leq (r_0-r)-(r_1-r_2)\leq C_{**}\,\sqrt\eps\,(s_0-s)\,(r_0-r_1)\,,
 \eeq

 where $C_{**}$ is a constant independent of $\eps$, $\rho_0$ and $\rho_1$.
 \ethm

 (A1) We shall use the following convention: $C_*$ denotes the maximum of the
 con-
 \hspace*{0.8cm} stants $C_*$ and $C_{**}$ from Theorems \ref{diffrieminv} and
 \ref{diffrieminv2}, respectively.\\

 In the following theorem $\beta$ and $\chi$ denote $S_1$ and
 $S_2$, respectively, while $o$ and $\pi$ denote $R_1$ and $R_2$,
 respectively. The prime is reserved for after interaction waves.
 (For example, the interaction of $S_2$ and $S_1$ which produces $S_1$
 and $S_2$ is denoted by $\chi+\beta\rightarrow
 \beta'+\chi'$.)

 \bthm\label{ocenatalasanew} If $0< \eps < \frac12$,
 then the following estimates are valid for the corresponding interactions:
 \begin{enumerate}
 \item $S_2$ and $S_1$ interaction:
 \begin{enumerate}
 \item $\chi+\beta\rightarrow\beta'+\chi'$\\
 $|\beta'|\leq |\beta|+C_*\sqrt\eps\,|\chi||\beta|,\;\;\;\; |\chi'|\leq
 |\chi|+C_*\sqrt\eps\,|\beta||\chi|$,\mbox{ or }\\

 there exist $\eta,\xi$ such that
 \item $\chi+\beta\rightarrow\beta'+\chi'$\\
 $0\leq |\beta'|=|\beta|-\xi,\;\;\;\; |\chi'|\leq
 |\chi|+C_*\sqrt\eps\,|\beta||\chi|+\eta$,\\
 where $0\leq \eta \leq g'_1(|\beta|,\rho_0)\xi<\xi$, or\\

 \item $\chi+\beta\rightarrow\beta'+\chi'$\\
 $0\leq |\chi'|=|\chi|-\xi,\;\;\;\; |\beta'|\leq
 |\beta|+C_*\sqrt\eps\,|\chi||\beta|+\eta$,\\
 where $0\leq \eta \leq g'_1(|\chi|,\rho_0)\xi<\xi$\;.
 \end{enumerate}

 \item $S_2$ and $R_1$ (or $R_2$ and $S_1$) interaction:
 \begin{enumerate}
 \item $\chi+o\rightarrow o'+\chi'$\\
 $|\chi'|=|\chi|,\;\;\;\; |o'|\leq |o|+C_*\sqrt\eps\,|\chi||o|$\;.
 \item $\pi + \beta\rightarrow \beta'+\pi'$\\
 $|\beta'|=|\beta|,\;\;\;\; |\pi'|\leq |\pi|+C_*\sqrt\eps\,|\beta||\pi|$\;.
 \end{enumerate}

 \item $S_2$ and $S_2$ (or $S_1$ and $S_1$) interaction:
 \begin{enumerate}
 \item $\chi_1+\chi_2\rightarrow o'+\chi':$\\
 $|\chi'|=|\chi_1|+|\chi_2|,\;\;\;\; |o'|\leq |\chi_1|+|\chi_2|$\;.
 \item $\beta_1+\beta_2\rightarrow \beta'+\pi':$\\
 $|\beta'|=|\beta_1|+|\beta_2|,\;\;\;\; |\pi'|\leq |\beta_1|+|\beta_2|$\;.
 \end{enumerate}

 \item $S_2$ and $R_2$ (or $R_1$ and $S_1$) interaction:
 \begin{enumerate}
 \item $1^\circ$ $\chi+\pi\rightarrow \beta'+\chi':$ there exist 1-shock
 $\beta_0$ and 2-shock $\chi_0$ such that\\
 $|\chi_0|=|\chi|-\xi,\;\;|\beta_0|=\eta$ and
 $\chi_0+\beta_0\rightarrow\beta'+\chi'$,\\
 where $0<\eta \leq g'_2(|\chi|,\rho_1)\xi<\xi$\;.\\
 $2^\circ$ $\chi+\pi\rightarrow \beta'+\pi'$: there exist $\eta,\xi$ such that\\
 $|\pi'|\leq |\pi|,\;\;\;\; |\beta'|=\eta<\xi=|\chi|$,\\
 where $0<\eta \leq g'_2(|\chi|,\rho_1)\xi<\xi$\;.
 \item $1^\circ$ $o+\beta\rightarrow \beta'+\chi':$ there exist 1-shock
 $\beta_0$ and 2-shock $\chi_0$ such that\\
 $|\beta_0|=|\beta|-\xi,\;\;|\chi_0|=\eta$ and $\chi_0+\beta_0\rightarrow\beta'+\chi'$,\\
 where $0<\eta \leq g'_2(|\beta|,\rho_2)\xi<\xi$\;.\\
 $2^\circ$ $o+\beta\rightarrow o'+\chi'$\\
 $|o'|\leq |o|,\;\;\;\; |\chi'|=\eta<\xi=|\beta|$,\\
 where $0<\eta \leq g'_1(|\beta|,\rho_0)\xi<\xi$\;.
 \end{enumerate}

 \item $R_2$ and $S_2$ (or $S_1$ and $R_1$) interaction:
 \begin{enumerate}
 \item $1^\circ$ $\pi+\chi\rightarrow \beta'+\chi':$ there exist $\eta,\xi$ such that\\
 $|\chi'|=|\chi|-\xi,\;\;\;\; |\beta'|=\eta$,\\
 where $0<\eta \leq g'_1(|\chi|,\rho_2)\xi<\xi$\;.\\
 $2^{\circ}$ $\pi+\chi\rightarrow \beta'+\pi'$\\
 $|\pi'|\leq |\pi|,\;\;\;\; |\beta'|=\eta<\xi=|\chi|$,\\
 where $0<\eta \leq g'_2(|\chi|,\rho_0)\xi<\xi$\;.
 \item $1^\circ$ $\beta+o\rightarrow \beta'+\chi':$ there exist $\eta,\xi$ such that\\
 $|\beta'|=|\beta|-\xi,\;\;\;\; |\chi'|=\eta$,\\
 where $0<\eta \leq g'_1(|\beta|,\rho_1)\xi<\xi$\;.\\
 $2^{\circ}$ $\beta+o\rightarrow o'+\chi'$\\
 $|o'|\leq |o|,\;\;\;\; |\chi'|=\eta<\xi=|\beta|$,\\
 where $0<\eta \leq g'_1(|\beta|,\rho_1)\xi<\xi$\;.
 \end{enumerate}

 \item $R_2$ and $R_1$ interaction:\\
 $\pi+o\rightarrow o'+\pi'$\\
 $|o'|=|o|,\;\;\;\; |\pi'|=|\pi|$.
 \end{enumerate}

 Here $C_*$ is a positive constant defined as in (A1).
 \ethm

 \pr This theorem can be proved using the same tools as in
 \cite{TNJS} and therefore will be omitted. The only differences are: the constant $C_*$ is now
 independent of $\eps,\beta,\chi, \rho_0, \rho_1$ and $\rho_2$,
 and we have $\sqrt\eps$ instead of $\eps$ in the estimates
 $(1)(a)-(1)(c)$, $(2)(a)$ and $(2)(b)$.
 \ep

 The main part of the paper is the interaction problem
 of delta shocks via pressure perturbation. Thus, one needs to control shock and rarefaction strengths as
 $\rho$ goes to infinity as $\eps\to 0$
 (more precisely, when $\rho$ is bounded by ${\rm const}/\eps$).
 Because of that, we give their estimates in $r-s$ plane
 based on Theorem \ref{diffrieminv} and Theorem \ref{diffrieminv2}.
 Let
 $(\rho_0,u_0)=(r_0,s_0)$ be connected with $(\rho,u)=(r,s)$ by a
 1-rarefaction (or 1-shock) wave, while $(\rho,u)=(r,s)$ be
 connected with $(\rho_1,u_1)=(r_1,s_1)$ by a 2-rarefaction (or 2-shock)
 wave. Then the strength of 1-rarefaction wave is
 \beq\label{strengthr1}
 r-r_0=\frac{2}{\sqrt\eps}(\rho_0^\eps-\rho^\eps),\;\; \rho<\rho_0\,,
 \eeq

 and the strength of 2-rarefaction wave is
 \beq\label{strengthr2}
 s_1-s=\frac{2}{\sqrt\eps}(\rho_1^\eps-\rho^\eps),\;\; \rho<\rho_1\,.
 \eeq

 The strength of 1-shock wave is estimated by
 \beq\label{strengths1}
 2\,\rho_0^\eps\,\sqrt\eps\,\ln\frac{\rho}{\rho_0}\leq r_0-r\leq
 2\,\frac{\sqrt\eps}{\sqrt{1+2\eps}}\,\left(\frac{\rho}{\rho_0}\right)^{\gamma/2}\cdot\rho_0^\eps,\;\;
 \rho>\rho_0\,,
 \eeq

 while, the strength of 2-shock wave is estimated by
 \beq\label{strengths2}
 2\,\rho_1^\eps\,\sqrt\eps\,\ln\frac{\rho}{\rho_1}\leq s-s_1\leq
 2\,\frac{\sqrt\eps}{\sqrt{1+2\eps}}\,\left(\frac{\rho}{\rho_1}\right)^{\gamma/2}\cdot\rho_1^\eps,\;\;
 \rho>\rho_1\,.
 \eeq

 Let us estimate the upper bound of the 1-shock wave given in (\ref{strengths1}).
 For the function $g_1$ from (\ref{g1fun}) we have $0\leq g'_1(\beta,\rho_0)<1$ and $0\leq
 g''_1(\beta,\rho_0)$, so
 $$
 \lim\limits_{|\beta|\to+\infty}g'_1(|\beta|,\rho_0)\leq 1.
 $$

 Let us consider two special cases needed for our investigation. The first case:
 $\rho>\rho_0$ and $\rho\sim 1/\eps$. We have that there exist
 constants $\tilde C$, $\bar{\bar C}$ and $\bar C$ independent of
 $\eps$ such that
 $$
 \ba{l}
 \ds\left(\frac{\rho}{\rho_0}\right)^{\gamma/2}\cdot\rho_0^\eps=
 \sqrt{\frac{\rho}{\rho_0}}\cdot\rho^\eps\leq\sqrt{\frac{\tilde C}{\eps}}\cdot\left(\frac{\bar{\bar C}}{\eps}\right)^\eps
 \leq\sqrt{\frac{1}{\eps}}\cdot \bar C,\;\;{\rm so }\\
 \\
 \ds
 2\,\frac{\sqrt\eps}{\sqrt{1+2\eps}}\,\left(\frac{\rho}{\rho_0}\right)^{\gamma/2}\cdot\rho_0^\eps\leq
 2\,\frac{\sqrt\eps}{\sqrt{1+2\eps}}\cdot\frac{\bar C}{\sqrt\eps}
 \leq{\rm{ const.}}
 \ea
 $$

 It follows that there exists a constant $C_2$, independent of $\eps$ and $\rho_0$, such
 that
 \beq\label{g1primocena}
 \sup g'_1(|\beta|,\rho_0):=C_2<1\;.
 \eeq

 Hence,
 \beq\label{g1primocena1}
 \frac{1-g'_1(|\beta|,\rho_0)}{g'_1(|\beta|,\rho_0)}\geq\frac{1-C_2}{C_2}=:C_3>0\,.
 \eeq

 The second case: $\rho>\rho_0$, $\rho\sim 1/\eps$ and $\rho_0\sim 1/\eps$.
 Then
 $$
 \sqrt{\frac{\rho}{\rho_0}}\cdot\rho^\eps\sim\rm{ const}\Rightarrow
 2\,\frac{\sqrt\eps}{\sqrt{1+2\eps}}\,\sqrt{\frac{\rho}{\rho_0}}\,\rho^\eps\sim
 {\mathcal O}(\sqrt\eps)\,.
 $$

 Again, $|\beta|\to\infty$ is impossible and
 (\ref{g1primocena1}) holds. In order
 to estimate the strength of $S_2$, we can use the same arguments
 to prove
 \beq\label{g2primocena}
 \sup g'_2(|\chi|,\rho_0)=:C_4<1,
 \eeq

 and
 \beq\label{g2primocena1}
 \ds
 \frac{1-g'_2(|\chi|,\rho_0)}{g'_2(|\chi|,\rho_0)}\geq\frac{1-C_4}{C_4}=:C_5>0\;.
 \eeq

 From now on, we shall put
 \beq\label{cnula}
 C_0=\min\{C_3,C_5\}\,.
 \eeq

 \section{Global interaction estimates}

 This section contains all the necessary assertions from
 \cite{FA} with several changes in constants. All changes are
 similar to those from the previous section.
 \bd\label{IJ krive novo} \cite{FA} A Lipschitz curve $J$ defined by $t=T(x)$, $x\in \mR$ is called an {\it
 I-curve}, if $|T'(x)|<1/\hat\lambda$. We denote $J_2>J_1$, if $T_1\neq T_2$
 and $T_2(x)\geq T_1(x), x\in \mR$. Denoting by $S_j(J)$ the set of $j$-shock
 waves crossing $J$ and $S(J)=S_1(J)+S_2(J)$, we define
 \beq\label{L}
 L^-(J)=\sum_{\alpha\in S(J)}|\alpha|,\;\;\;\;Q(J)=\sum_{\beta\in S_1(J),\chi\in
 S_2(J),\;\;\beta,\chi\mbox{\tiny{ approach}}}|\beta||\chi|\;.
 \eeq
 \ed

 Set $F(J)=L^-(J)+\tilde K\cdot Q(J)$, where $\tilde K:=4C_*\sqrt\eps$.
 A space-like line lying between the initial
 line and the first interaction point is denoted with $O$.

 \blem\label{q0}
 \beq\label{q0ocena}
 Q(O)\leq L^-(O)^2\;.
 \eeq

 \elem

 \pr The proof follows straightforward from Definition \ref{IJ krive novo}.
 \ep

 \blem\label{f0}
 Assuming $\;4\,C_*\,\sqrt\eps\,L^-(O)\leq 1$, we have
 \beq\label{f0ocena}
 F(O)\leq 2L^-(O)\;.
 \eeq
 \elem

 \pr
 $$
 \ba{rl}
 F(O)=&L^-(O)+\tilde K\;Q(O)\leq L^-(O)+\tilde K\;L^-(O)^2 \mbox{ (by (\ref{q0ocena}))}\\
 =&L^-(O)(1+\tilde KL^-(O))=L^-(O)(1+4\;C_*\,\sqrt\eps\,L^-(O))\\
 \leq& L^-(O) (1+1)=2L^-(O)\;.
 \ea
 $$
 \ep

 As in \cite{AB}, consider a interval ${\mathcal J}\subset \mR$
 and a map $a:{\mathcal J}\to\mR^n$. The {\it total variation} (TV) of $a$ is
 then defined as
 $$
 TV (a):=\sup\left\{\sum_{j=1}^N|a(x_j)-a(x_{j-1})|\right\},
 $$

 where the supremum is taken over all $N\geq 1$ and all
 $(N+1)$-tuples of points $x_j\in {\mathcal J}$ such that
 $x_0<x_1<\cdots<x_N$. Now, we give a new estimate for $L^+(O)$.
 Here, $L^+(O)$ denotes the sum of the
 rarefaction waves strengths which cross the line $O$.
 \blem\label{l0tv} We have
 \beq\label{l0tvocena}
 L^-(O)\leq TV(r_0(x),s_0(x))\;\;\mbox{ and }\;\;L^+(O)\leq
 TV(r_0(x),s_0(x))\;.
 \eeq
 \elem

 The estimates in previous Lemma can easily be verified.
 The uniform bounds of $F(J)$ follows from the following theorem.
 \bthm\label{thm5} If $\ds\; C_*\,\sqrt\eps\,F(O)\leq\min\left\{\frac{1}{2},\frac{C_0}{4}\right\}$, then
 $F(J_2)\leq F(J_1)$ for $J_2>J_1$. Particulary, $L^-(J)\leq
 F(O)$.
 \ethm

 \pr This theorem can be proved in the same way as Lemma 5 from \cite{TNJS} and
 hence the proof will be omitted.
 One has just to substitute a constant $K$ from the original proof with
 the determined value $\tilde K$ here.
 \ep

 \blem\label{lemaAsakura} Assume that $\tilde K L^-(O)\leq 1$ and that
 \beq\label{mainas}
 \sqrt{\gamma-1}\;TV\,(r_0(x),s_0(x))\leq
 \frac1{C_*}\cdot\min\left\{\frac{\sqrt2}{4},\frac{\sqrt2}{8}\,C_0\right\}\,.
 \eeq

 Then $\ds \tilde K F(O)\leq \min\left\{2,C_0\right\}$.
 \elem

 \pr Using Lemma \ref{f0} and \ref{l0tv} we have
 $$
 \frac{\sqrt2}{2}\sqrt\eps\,F(O)\leq \sqrt2\sqrt\eps\,L^-(O)\leq \sqrt{\gamma-1}\;TV\,(r_0(x),s_0(x))\leq
 \frac1{C_*}\cdot\min\left\{\frac{\sqrt2}{4},\frac{\sqrt2}{8}\,C_0\right\}\,.
 $$

 Multiplying it with $8\,C_*/\sqrt{2}$, one gets
 $$
 4\,C_*\sqrt\eps\,F(O)=\tilde KF(O)\leq \min\left\{2,C_0\right\}\,.
 $$

 which proves the claim.
 \ep

 The right hand side of (\ref{mainas}) does not depend
 on $\eps$, and then one can say that $TV\,(r_0(x),s_0(x))$ may be arbitrarily large since we can always
 choose $\eps$ small enough in order to fulfill (\ref{mainas}) with $\gamma=1+2\eps$.
 Then we can apply wave front tracking procedure from \cite{FA} for each such $\eps$,
 and obtain a sequence of step functions converging to the
 entropic solution. One
 only needs to replace $C\eps F(O)$ and
 $\frac{1-\delta}{\delta}$ from \cite{FA} with $C_*\sqrt\eps F(O)$ and
 $C_0$, respectively.

 \section{Approximate delta shock solutions to pressureless gas dynamics}

 Our main task is to solve delta shock interaction problem for
 pressureless gas dynamics model. Accordingly, we will introduce a solution
 concept from \cite{MN2009p} (somewhat simplified) and check consistency of
 theoretical and numerical wave front tracking results by letting $\eps\to 0$.

 \subsection{Basic notions}
 In this section we shall use the notions and assertions from \cite{MN2009p}.
 It contains results for a $3\times 3$ system with energy
 conservation law added, but all the results can also be applied to system (\ref{PGD}), too.
 Let us start with the basic definitions. Vector valued function of the form
 \begin{equation} \label{sdw}
 U_{\varepsilon}(x,t)=\begin{cases} U_{0},
 & x<c(t)-a_{\varepsilon}(t) \\
 U_{1,\varepsilon}(t), & c(t)-a_{\varepsilon}(t)<x<c(t) \\
 U_{2,\varepsilon}(t), & c(t)<x<c(t)+b_{\varepsilon}(t) \\
 U_{1}, & x>c(t)+b_{\varepsilon}(t)\end{cases}.
 \end{equation}

 is called {\it weighted shadow wave} (weighted SDW, for short).
 Here, $U:=(\rho,u)$. The functions $a_{\varepsilon}$, $b_{\varepsilon}$ are
 continuous functions satisfying
 $a_{\varepsilon}(0)=x_{1,\varepsilon}$ and
 $b_{\varepsilon}(0)=x_{2,\varepsilon}$.
 The SDW is {\it constant} if $U_{1,\varepsilon}$ and
 $U_{2,\varepsilon}$ are just constants. If, in addition,
 $x_{1,\varepsilon}=x_{2,\varepsilon} =0$, then the wave is called
 {\it simple}.

 The value
 \begin{equation*}
 \sigma_{\varepsilon}(t):=a_{\varepsilon}(t)U_{1,\varepsilon}(t)
 +b_{\varepsilon}(t)U_{2,\varepsilon}(t)
 \end{equation*}
 is called the {\it strength} and $c'(t)$ is called the {\it speed}
 of the shadow wave. We assume that $\lim_{\varepsilon \rightarrow
 0} \sigma_{\varepsilon}(t)=\sigma(t)\in \mathbb{R}^{n}$ exists for
 every $t\geq 0$ and
 \begin{equation*}
 \lim_{\varepsilon \rightarrow 0}\int
 U_{\varepsilon}(x,t)\phi(x,t)\,dx\,dt
 =\langle U_0+(U_1-U_0)\,\theta(x-c(t))+\sigma(t)\,\delta(x-c(t)),\phi(x,t)\rangle,
 \end{equation*}

 for $t\geq 0$, where $\theta$ is a Heaviside function. The SDW {\it central line} is given by $x=c(t)$, while
 $x=c(t)-a_{\varepsilon}(t)$ and $x=c(t)+b_{\varepsilon}(t)$ are
 called the {\it external SDW lines}. The values
 $x_{1,\varepsilon}$ and $x_{2,\varepsilon}$ are called the {\it
 shifts}, while $U_{1,\varepsilon}(t)$ and $U_{2,\varepsilon}(t)$
 are called the {\it intermediate states} of a given SDW.

 Let $i\in\{1,2,\dots,n\}$. We assume $
 \|U_{\varepsilon}^{i}\|_{L^{\infty}}= {\mathcal
 O}(\varepsilon^{-1}),
 $ if $f$ and $g$ have at most
 a linear growth with respect to $i$-th component, or otherwise
 $
 \|U_{\varepsilon}^{i}\|_{L^{\infty}}= {o}(\varepsilon^{-1})$.
 The components of the first kind are called {\it major} ones, while
 the ones of the second kind are called {\it minor} ones.

 A {\it delta shock} is a SDW associated with a $\delta$
 distribution with all minor components having finite limits as
 $\varepsilon \rightarrow 0$.
 \medskip

 The following lemma is the base of all calculations involving
 SDWs.

 \blem \label{tezinska} Let $f,g\in {\mathcal C}(\Omega:{\mathbb
 R}^{n})$ and $U:{\mathbb R}_{+}^{2}\rightarrow \Omega\subset
 {\mathbb R}^{n}$ be a piecewise constant function for every $t\geq
 0$. Let us also suppose that $f$ and $g$ satisfy
 \begin{equation}\label{gass*}
 \max_{i=1,2}\{\|f(U_{i,\varepsilon})\|_{L^{\infty}},
 \|g(U_{i,\varepsilon})\|_{L^{\infty}}\} ={\mathcal
 O}(\varepsilon^{-1}).
 \end{equation}
 Then
 \begin{equation}\label{gosnrelft}
 \begin{split}
 \langle \partial_{t}f(U_{\varepsilon}),\phi \rangle \approx &
 \int_{0}^{\infty} \lim_{\varepsilon \rightarrow 0} {d \over
 dt}\Big(a_{\varepsilon}(t)f(U_{1,\varepsilon}(t))
 +b_{\varepsilon}(t)f(U_{2,\varepsilon}(t))\Big) \, \phi(c(t),t) \, dt \\
 & -\int_{0}^{\infty} c'(t) \Big(f(U_{1})-f(U_{0})\Big)\,
 \phi(c(t),t) \, dt \\
 & +\int_{0}^{\infty}\lim_{\varepsilon \rightarrow 0}
 c'(t)\Big(a_{\varepsilon}(t)f(U_{1,\varepsilon}(t))
 +b_{\varepsilon}(t)f(U_{2,\varepsilon}(t))\Big) \,
 \partial_{x}\phi(c(t),t) \, dt
 \end{split}
 \end{equation}
 and
 \begin{equation}\label{gosnrelgt}
 \begin{split}
 \langle \partial_{x}g(U_{\varepsilon}),\phi \rangle \approx
 & \int_{0}^{\infty}\Big(g(U_{1})-g(U_{0})\Big)\, \phi(c(t),t) \, dt \\
 & -\int_{0}^{\infty}\lim_{\varepsilon \rightarrow 0}
 \Big((a_{\varepsilon}(t)g(U_{1,\varepsilon}(t))
 +(b_{\varepsilon}(t)g(U_{2,\varepsilon}(t))\Big)\,
 \partial_{x}\phi(c(t),t)\, dt.
 \end{split}
 \end{equation}
 \elem

 \subsection{Entropy conditions}

 Let $\eta(U)$ be a {\it semi-convex} entropy function for
 (\ref{PGD}), with entropy-flux function $q(U)$. We shall use
 entropy condition in the following form. A weak or approximate
 solution $U_{\varepsilon}=(\rho_{\varepsilon},u_{\varepsilon})$ to
 system (\ref{PGD}) with initial data
 $U|_{t=0}=U_{0,\varepsilon}$ is {\it admissible} provided that for every
 $T>0$ we have
 \begin{equation} \label{ent-osn}
 \underline{\lim}_{\varepsilon \rightarrow 0} \int_{\mathbb
 R}\int_{0}^{T} \eta(U_{\varepsilon}) \partial_{t}\phi
 +q(U_{\varepsilon}) \partial_{x}\phi \, dt \, dx +\int_{\mathbb R}
 \eta(U_{0,\varepsilon}(x,0))\phi(x,0) \, dx \geq 0,
 \end{equation}
 for all non-negative test functions $\phi\in
 C_{0}^{\infty}(\mathbb{R}\times (-\infty,T))$.

 Using Lemma \ref{tezinska} with $f$ substituted by $\eta$ and $g$
 by $q$ and the fact that the delta function is a non-negative
 distribution, the first condition for SDW $U_{\varepsilon}$ from
 (\ref{sdw}) to be admissible is given by
 \begin{equation} \label{E1}
 \begin{split}
 -c'(t)(\eta(U_{1})-\eta(U_{0}))+(q(U_{1})-q(U_{0})) & \\
 +\lim_{\varepsilon \rightarrow 0}{d \over dt}
 (\eta(U_{1,\varepsilon}(t))a_{\varepsilon}
 +\eta(U_{2,\varepsilon}(t))b_{\varepsilon}) & \leq 0.
 \end{split}
 \end{equation}
 The derivative of delta function changes the sign, so
 $U_{\varepsilon}$ has to satisfy
 \begin{equation} \label{E2}
 \begin{split}
 \lim_{\varepsilon \rightarrow 0}c'(t)
 (\eta(U_{1,\varepsilon}(t))a_{\varepsilon}
 +\eta(U_{2,\varepsilon}(t))b_{\varepsilon})& \\
 -q(U_{1,\varepsilon}(t))a_{\varepsilon}(t)
 -q(U_{2,\varepsilon}(t))b_{\varepsilon}(t) & = 0
 \end{split}
 \end{equation}
 in addition.

 These conditions are much simpler in the case of simple SDW
 when $U_{0}$, $U_{1}$, $U_{1,\varepsilon}$ and
 $U_{2,\varepsilon}$ are constants:
 \begin{equation}\label{entrouslovi1}
 \overline{\lim}_{\varepsilon \rightarrow 0}
 -c(\eta(U_{1})-\eta(U_{0}))+
 a_{\varepsilon}\eta(U_{1,\varepsilon})+b_{\varepsilon}\eta(U_{2,\varepsilon})
 +q(U_{1})-q(U_{0})\leq  0
 \end{equation}
 and
 \begin{equation} \label{entrouslovi2}
 \lim_{\varepsilon \rightarrow 0}
 -c(a_{\varepsilon}\eta(U_{1,\varepsilon})
 +b_{\varepsilon}\eta(U_{2,\varepsilon}))
 +a_{\varepsilon}q(U_{1,\varepsilon})+b_{\varepsilon}q(U_{2,\varepsilon})
 = 0.
 \end{equation}

In most of the papers with delta or singular shock solution, the
authors use overcompressibility as the admissibility condition. A
wave is called the {\it overcompressive} one if all
characteristics from both sides of the SDW line run into a shock
curve, i.e.
\begin{equation*}
\lambda_{i}(U_{0})\geq c'(t) \geq \lambda_{i}(U_{1}),\;
i=1,\ldots,n,
\end{equation*}
where $c$ is a shock speed and $x=\lambda_{i}(U)t$, $i=1,\ldots,n$
are the characteristics of the system. One will see that these
notations coincide with our model case.

The entropy condition is connected with the problem of uniqueness
for a weak solution of the conservation law system. We give a
definition of weak (distributional) uniqueness and some results
about it afterward.

\bd \label{d-uniqueness} An SDW solution is called {\it weakly
unique} if its distributional image is unique. More precisely, a
speed $c$ of the wave has to be unique as well as the limit
\begin{equation*}
\lim_{\varepsilon \rightarrow 0} a_{\varepsilon}U_{1,\varepsilon}
+b_{\varepsilon}U_{2,\varepsilon}.
\end{equation*}
Let $i\in \{1,\ldots,n\}$. If a limit $\displaystyle
\lim_{\varepsilon \rightarrow 0}
a_{\varepsilon}U_{1,\varepsilon}^{i}+b_{\varepsilon}U_{2,\varepsilon}^{i}$
is unique, then we say that the $i$-th component is unique. \ed

Note that all minor components of $U_{\varepsilon}$ are unique by
default.

\subsection{Entropy solutions to Riemann problem for
pressureless gas dynamics model}

The proof for the following theorem in the case of $3\times 3$ PGD
model is given in \cite{MN2009p}. Its restriction to a $2\times 2$
system is straightforward and therefore not discussed here.

\bthm \label{tPGD} Suppose that $u_{0}>u_{1}$. Then there exists a
unique shadow wave solution of the form (\ref{sdw}) to the Riemann
problem (\ref{PGD}, \ref{idPGD}) satisfying the entropy inequality
(\ref{ent-osn}) with $\eta$ and $q$ as defined above.

Moreover, the validity of (\ref{ent-osn}) for all semi-convex
entropies $\eta$ are equivalent to the overcompressibility of the
shadow wave. \ethm

Our aim is to show the structure of a solution in order to be able
to compare it with a numerical approximation described above. For
our purposes it is safe to take
$a_{\varepsilon}=b_{\varepsilon}=\varepsilon$ in the sequel. In
the proof of Theorem \ref{tPGD} we showed that a SDW solution
(\ref{sdw}) (with $U=(\rho,u)$) to (\ref{PGD}) and initial data
(\ref{idPGD}), with $u_{0}>u_{1}$, had to satisfy
\begin{equation*}
\begin{split}
c=u_s=\lim_{\varepsilon \rightarrow 0}u_{\varepsilon} & \equiv
\frac{[\rho u]-[u]\sqrt{\rho_{0}\rho_{1}}}{[\rho]}
\; (u_s\text{ does not depend on } \varepsilon)\\
\lim_{\varepsilon \rightarrow 0}\varepsilon \rho_{\varepsilon} &=
c[\rho]-[\rho u]=(u_{0}-u_{1})\sqrt{\rho_{0}\rho_{1}},
\end{split}
\end{equation*}

if $\rho_0\neq \rho_1$, and $c=u_s=(u_0+u_1)/2$, if
$\rho_0=\rho_1$. That defines a weakly unique SDW solution to the
problem.

\subsection{Two SDWs interaction}\label{sdwsdwint}

The main advantage of using weighted SDWs (intermediate states
vary with $t$ in addition) is for solving SDW interaction problem.
Then we can proceed with the main part of the paper by showing
numerically that such a solution can be viewed as a limit of gas
dynamics model with a vanishing pressure as perturbation. Note
that verification of delta shock existence has already been
obtained in \cite{GQCHL} (see \cite{DMMN} for a somewhat general
model).

Suppose that two SDWs interact in a point $(X,T)$. The superscript
$1$ is used for data in the left wave while the superscript $2$ is
used for the right one. The first SDW connects the states
$\displaystyle U_{0}=(\rho_{0},u_{0})$ with $\displaystyle
U_{1}=(\rho_{1},u_{1})$, while the second one connects the states
$\displaystyle U_{1}=(\rho_{1},u_{1})$ with $\displaystyle
U_{2}=(\rho_{2},u_{2})$.

Again, the following theorem has been proved in \cite{MN2009p} for
the extended PGD system, and the proof can easily be adopted for
the present one (\ref{PGD}).

\bthm \label{t-varSDW} The result of two SDW interactions for the
pressureless system (\ref{PGD}) is a weakly unique single entropic
weighted SDW. \ethm

We use the following notation: $[x]_1:=x_1-x_0$, $[x]_2:=x_2-x_1$
and $[x]:=x_2-x_0$. The weighted SDW solution from the above
theorem satisfies the following: The speed is given by
$c'(t)=u_s(t):=\lim_{\varepsilon \rightarrow
0}u_{\varepsilon}(t)$, while $u_s(t)$ and
$\xi(t):=\lim_{\varepsilon \rightarrow 0}\varepsilon
\rho_{\varepsilon}(t)$ satisfies the following ODEs system
\begin{equation} \label{sdwsdw}
\begin{split}
\xi'(t) &= u_s(t)[\rho]-[\rho u] \\
(\xi(t)u_s(t))'&=u_s(t)[\rho u]-[\rho u^{2}]
\end{split}
\end{equation}
with the initial data
\beq\label{sdwsdwid}
\begin{split}
\xi(T)= & (\xi^{1}+\xi^{2})T=(-[u]_1\sqrt{\rho_{0}\rho_{1}}
-[u]_2\sqrt{\rho_{1}\rho_{2}})T, \\
\xi(T)u_{s}(T)= &(c^{1}\xi^{1}+c^{2}\xi^{2})T =\Big(-\frac{[\rho
u]_1-[u]_1\sqrt{\rho_{0}\rho_{1}}}
{[\rho]_1}\cdot [u]_1\sqrt{\rho_{0}\rho_{1}}\\
& - \frac{[\rho u]_2-[u]_2\sqrt{\rho_{1}\rho_{2}}}
{[\rho]_2}\cdot[u]_2\sqrt{\rho_{1}\rho_{2}}\Big)T.
\end{split}
\eeq

Here are some facts regarding the solution $(\xi(t),u_s(t))$,
$t\geq T$ to the above initial data problem (see \cite{MN2009p}):
\begin{enumerate}
\item $\xi(t)$, for $t>T$, is an increasing function when exists.
The initial data $\xi(T)>0$ and $\xi(t)$ is always positive
function for $t>T$ (when exists), since $u_0>u_1>u_2$.
\item From the system (\ref{sdwsdw}) we have
$$
u_s'(t)=-\frac{1}{\xi(t)}([\rho]u_s^2(t)-2[\rho u]u_s(t)+[\rho
u^2]).
$$

The value $-1/\xi(t)$ is now always negative for $t>T$. The roots
of the right-hand side of the above ODE are denoted as $A_1<A_2$.
Then, for $[\rho]\neq 0$,
$$
A_{1,2}=\frac{[\rho u]\pm |u_0-u_2|\sqrt{\rho_0\rho_2}}{[\rho]}\,.
$$
Assume that $[\rho]>0$. If $u_s(t)\in (A_1,A_2)$, then $u_s(t)$
increases, and if $u_s(t)\in (-\infty,A_1)\cup(A_2,+\infty)$, then
$u_s(t)$ decreases. The opposite holds if $[\rho]<0$. There are
two possible cases:
\begin{itemize}
\item If $\rho_0>\rho_2$, then $u_2\leq A_1\leq u_0\leq A_2$. If
$u_s(T)\in(u_2,A_1)$, then $u_s(t)$ increases for $t>T$ but stays
bellow $A_1$. If $u_s(T)\in(A_1,u_0)$, then $u_s(t)$ decreases for
$t>T$ but stays above $A_1$.
\item If $\rho_2>\rho_0$, then $A_1\leq u_2\leq A_2\leq u_0$. Again, if
$u_s(T)\in(u_2,A_2)$, then $u_s(t)$ increases for $t>T$ but stays
bellow $A_2$. If $u_s(T)\in(A_2,u_0)$, then $u_s(t)$ decreases for
$t>T$ but stays above $A_2$.
\end{itemize}
This implies $u_0\geq u_s(t)\geq u_2$ (the SDW is
overcompressive). Also, one will see that numerical examples
resemble these asymptotic properties of $u_s(t)$ as $t\to\infty$.
\end{enumerate}

 \section{Numerical results}

 In this section one can find numerical results which show a
 consistency of theoretical (in the sense of SDWs) and numerical results.
 Consider system (\ref{modsys}) with the initial data
 \beq\label{initconnum}
 (\rho,u)\left|_{t=0}\right.=\left\{
 \ba{ll}
 (\rho_0,u_0), & x<a_1\\
 (\rho_1,u_1), & a_1<x<a_2\\
 (\rho_2,u_2), & x>a_2
 \ea
 \right.
 \eeq

 where $a_1<a_2$, $u_0>u_1>u_2$. Then (see \cite{GQCHL}), for $\eps$ small
 enough, there exist $(\rho_{1,\eps},u_{1,\eps})\in \mR_+\times
 \mR$ and $(\rho_{2,\eps},u_{2,\eps})\in \mR_+\times
 \mR$, so that:
 \begin{itemize}
 \item $(\rho_0,u_0)$ is connected with
 $(\rho_{1,\eps},u_{1,\eps})$ by an 1-shock, and
 $(\rho_{1,\eps},u_{1,\eps})$ is connected with $(\rho_1,u_1)$ by
 a 2-shock,
 \item $(\rho_1,u_1)$ is connected with
 $(\rho_{2,\eps},u_{2,\eps})$ by an 1-shock, while
 $(\rho_{2,\eps},u_{2,\eps})$ is connected with $(\rho_2,u_2)$ by
 a 2-shock.
 \end{itemize}

 A numerical solution is obtained by wave front tracking algorithm described in
 \cite{FA}. In order to verify two delta shocks interaction, we shall consider two
 cases.\\

 {\bf Case A.} Suppose that
 $(\rho_0,u_0)$ is connected with $(\rho_1,u_1)$ by a single delta shock and
 $(\rho_1,u_1)$ is connected with $(\rho_2,u_2)$ by a single
 delta shock, too. Assume that $(\rho_0,u_0)$ can be connected with $(\rho_2,u_2)$ by a single
 delta shock (so-called simple SDW, see \cite{MN2009p}).
 The resulting SDW has a constant speed as a consequence. That can be
 done by choosing a special value for $\rho_2$ provided that
 $\rho_0,u_0,\rho_1,u_1$ and $u_2$ are already given.\\
 {\bf Case B.} We choose arbitrarily $\rho_2$, i.e. the resulting SDW has a variable speed
 (a central SDW curve is no longer a line). The numerical results
 are given in Tables \ref{ds1}, \ref{ds2} and \ref{ds3}.

 \begin{center}
 \begin{table}
 \caption{Parameter description}\label{pardeskr}
 \begin{tabular}[htbp]{|l|l|}\hline
 Parameter & Description\\\hline\hline
 $\kappa$ & Adiabatic constant defined in (\ref{definP}).\\\hline
 $\rho_\eps$ & First component of the intermediate state of the solution for (\ref{modsys}, \ref{initconnum}).\\\hline
 $u_\eps$ & Second component of the intermediate state of the solution for (\ref{modsys}, \ref{initconnum}).\\\hline
 $c_1$ & Speed of the first left shock.\\\hline
 $c_2$ & Speed of the last right shock.\\\hline
 $|Eq_1|$ & Left hand side of the integral on the first equation in (\ref{modsys}).\\\hline
 $|Eq_2|$ & Left hand side of the integral on the second equation in (\ref{modsys}).\\\hline
 \end{tabular}
 \end{table}
 \end{center}

 \subsection{Case A}
 \bex\label{ex1}
 Let $a_1=0$, $a_2=2$, $(\rho_0,u_0)=(1,1)$, $(\rho_1,u_1)=(1.2,0.8)$ and $u_2=0.7$.
 Now, for $\rho_2=1.14286$, there exists a single simple SDW as a solution to the interaction
 problem.
 \eex

 \begin{table}[htbp]
 \caption{$(\rho_0,u_0)=(1,1),(\rho_1,u_1)=(1.2,0.8),(\rho_2,u_2)=(1.14286,0.7)$}\label{ds1}
 \vspace*{-0.2cm}
 $$
 \begin{array}{c|c|c|c|c|c|c|c|c|c|}
 \gamma & \kappa & \eps  & \rho_{\eps} &  u_\eps &   c_1   &   c_2   &  |Eq_1|         & |Eq_2|         \\\hline
 2      & 0.5    & 0.5   & 1.29979     & 0.80062 & 0.13553 & 1.53337 &  2\cdot 10^{-5} & 3\cdot 10^{-5} \\\hline
 1.2    & 0.29   & 0.1   & 1.68612     & 0.82806 & 0.57746 & 1.09746 &  1\cdot 10^{-4} & 2\cdot 10^{-4} \\\hline
 1.02   & 0.099  & 0.01  & 4.22718     & 0.84163 & 0.79256 & 0.89411 &  2\cdot 10^{-3} & 3\cdot 10^{-3} \\\hline
 1.01   & 0.071  & 0.005 & 6.71337     & 0.84311 & 0.81565 & 0.87247 &  8\cdot 10^{-5} & 2\cdot 10^{-3} \\\hline
 1.006  & 0.055  & 0.003 & 9.95729     & 0.84379 & 0.82636 & 0.86254 &  1\cdot 10^{-2} & 6\cdot 10^{-3} \\\hline
 \end{array}
 $$
 \end{table}

\begin{figure}[htbp]
 \centerline{
 \epsfig{figure=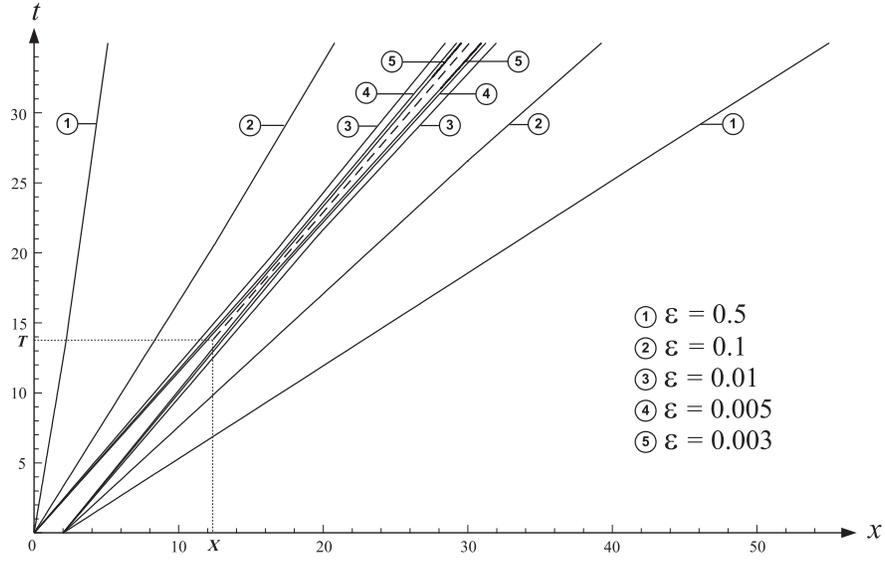,width=12cm
 }} \epsfxsize=1\textwidth
 \caption{Phase $x-t$ plane, Case A, Example \ref{ex1}.\label{test}}
 \end{figure}

 After interaction, the speed of the resulting wave is $c_\delta=0.844994$. Two SDWs will interact in a point
 $(X,T)=(12.365,13.8088)$ with such data.
 Now, we are going to explain Figures \ref{test}, \ref{test2} and \ref{test3}
 which are illustrations of appropriate numerical results.
 For each $\eps$ we have two piecewise linear half-lines. The left one originates from the point $(x,t)=(a_1,0)$,
 while the right one originates from the point $(x,t)=(a_2,0)$.
 The $i$-th linear segment of these half-lines can be written in the form
 $x=c_{i,j}\,(t-t_i)+x_i$, $i\geq 1$, $j=1,2$, $x_i\leq x\leq x_{i+1}$, $t_i\leq t\leq t_{i+1}$,
 where $c_{i,1}$ stands for the speed of the first ($S_1$) wave on the
 left hand side in phase plane, while
 $c_{i,2}$ stands for the speed of the last ($S_2$) wave on the
 left hand side in phase plane at each $i$-th
 segment. Interactions of the waves occur at the points
 $(x_i,t_i),\;i\geq 1$. After two delta shock interaction, the resulting delta shock central line in Figure
 \ref{test} (dashed line) starts from $(X,T)$ and it is
 calculated explicitly from system (\ref{PGD}).

 \subsection{Case B}

\begin{figure}[htbp]
 \centerline{
 \epsfig{figure=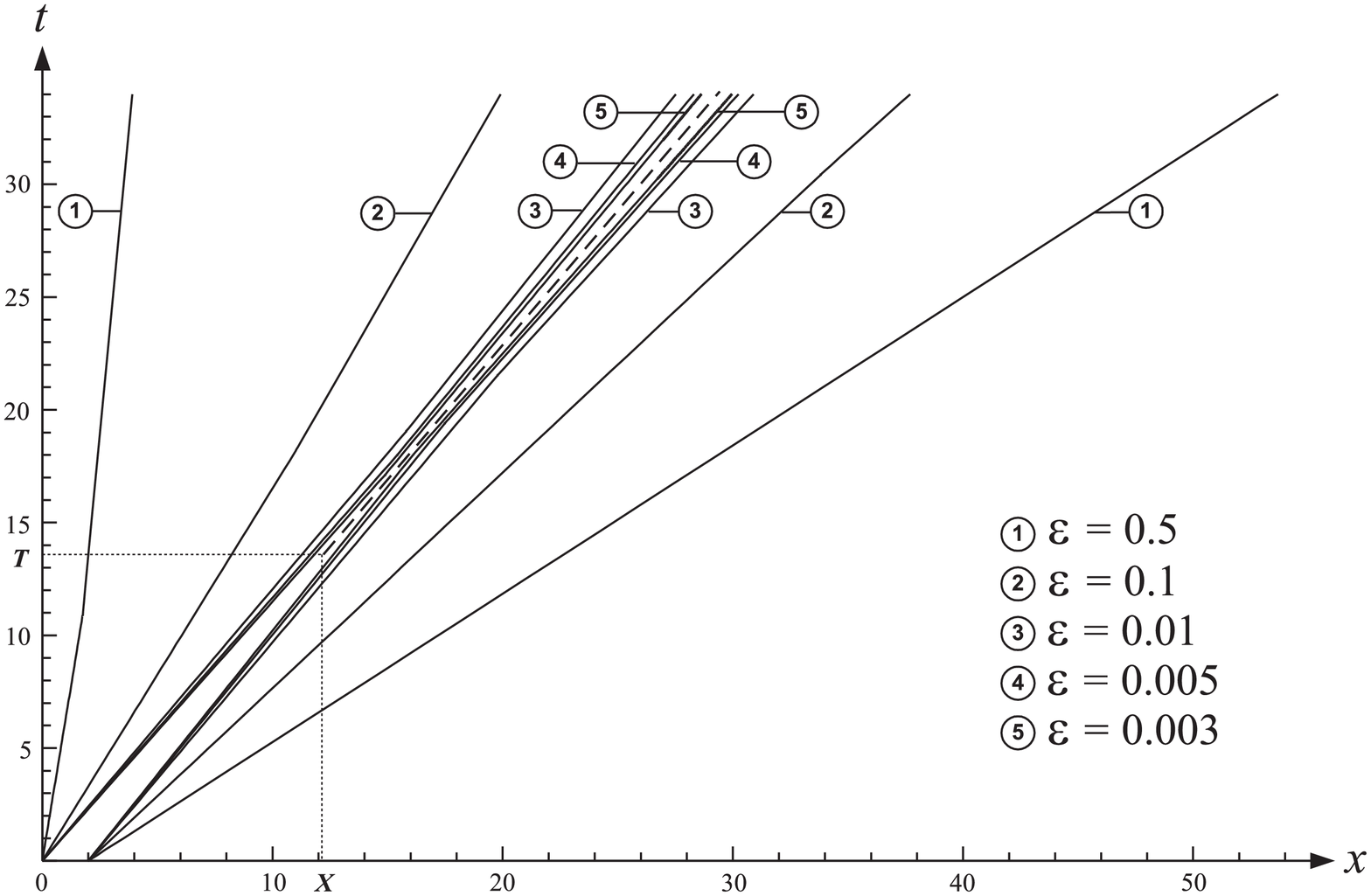,width=12cm
 }} \epsfxsize=1\textwidth
 \caption{Phase $x-t$ plane, Case B, Example \ref{ex2}.\label{test2}}
 \end{figure}

\begin{figure}[htbp]
 \centerline{
 \epsfig{figure=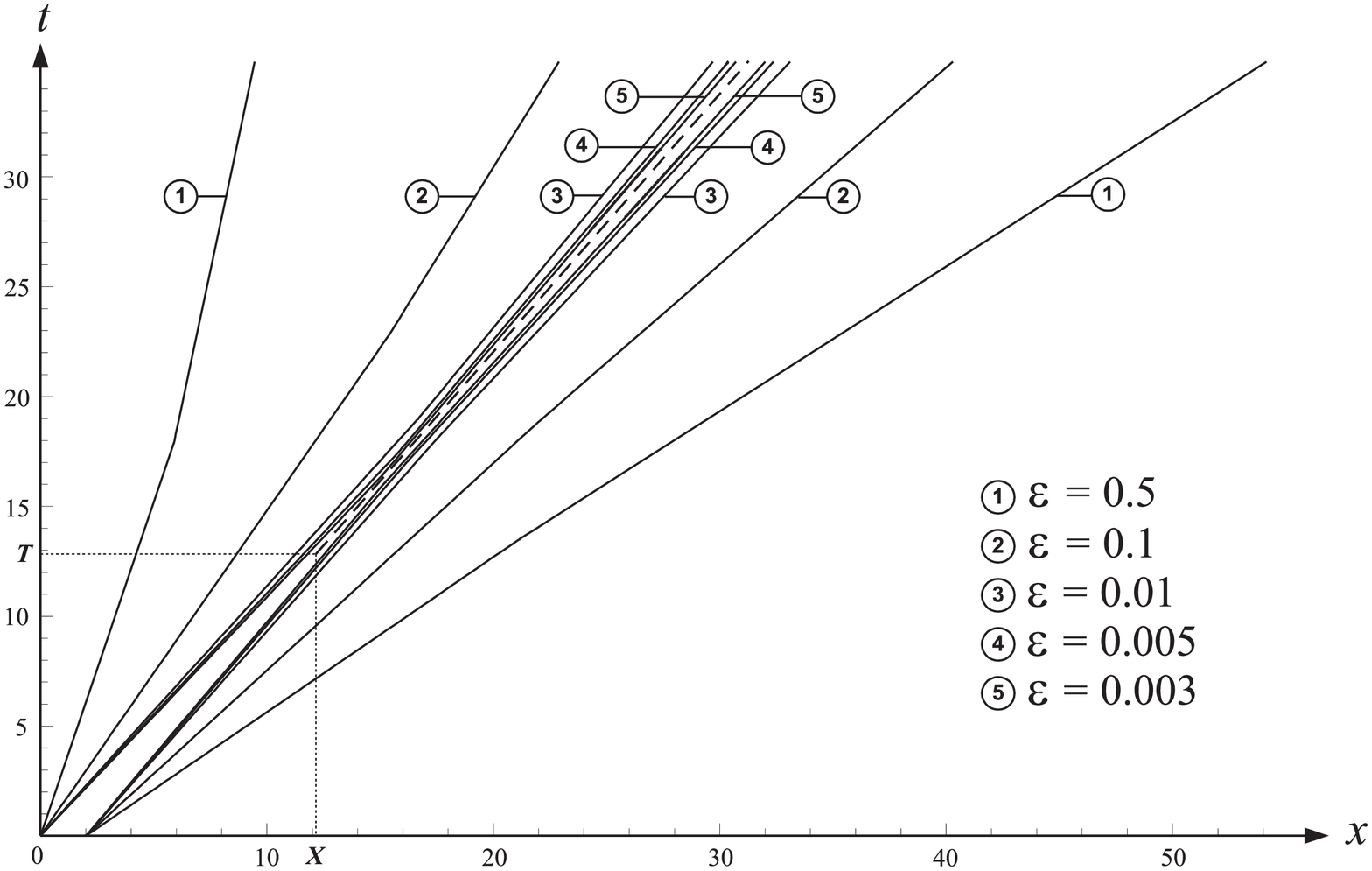,width=12cm
 }} \epsfxsize=1\textwidth
 \caption{Phase $x-t$ plane, Case B, Example \ref{ex3}.\label{test3}}
 \end{figure}

 \bex\label{ex2}
 Let $a_1=0$, $a_2=2$, $(\rho_0,u_0)=(1,1)$, $(\rho_1,u_1)=(1.2,0.8)$ and
 $(\rho_2,u_2)=(1.3,0.7)$. Two SDWs will interact in a point
 $(X,T)=(12.2291,13.657)$ with such data. After two delta shock interaction, the resulting delta shock central lines
 in Figures \ref{test2} and \ref{test3} (dashed lines) start from $(X,T)$, too.
 Here,
 $$
 x(t)=\int\limits_T^t u_s(p)\,dp+X,\;t\geq T,
 $$
 while $u_s(t)$ represents the second component of the solution $(\xi(t),u_s(t))$ of system
 (\ref{sdwsdw})  with initial conditions (\ref{sdwsdwid}).
 \eex

 \begin{table}[htbp]
 \caption{$\;\;(\rho_0,u_0)=(1,1),(\rho_1,u_1)=(1.2,0.8),(\rho_2,u_2)=(1.3,0.7)$}\label{ds2}
 \vspace*{-0.3cm}
 $$
 \begin{array}{c|c|c|c|c|c|c|c|c|c|}
 \gamma & \kappa & \eps  & \rho_\eps &  u_\eps &   c_1   &   c_2   &  |Eq_1|         & |Eq_2| \\\hline
 2      & 0.5    & 0.5   & 1.38131   & 0.74968 & 0.09312 & 1.54395 &  3\cdot 10^{-4} & 5\cdot 10^{-4} \\\hline
 1.2    & 0.29   & 0.1   & 1.79287   & 0.80661 & 0.56268 & 1.08778 &  5\cdot 10^{-4} & 6\cdot 10^{-4} \\\hline
 1.02   & 0.099  & 0.01  & 4.49667   & 0.83356 & 0.78596 & 0.88787 &  3\cdot 10^{-3} & 3\cdot 10^{-3} \\\hline
 1.01   & 0.071  & 0.005 & 7.14038   & 0.83646 & 0.80983 & 0.86684 &  1\cdot 10^{-3} & 3\cdot 10^{-3} \\\hline
 1.006  & 0.055  & 0.003 & 10.5884   & 0.83782 & 0.82091 & 0.85711 &  3\cdot 10^{-2} & 2\cdot 10^{-2} \\\hline
 \end{array}
 $$
 \end{table}

 \bex\label{ex3}
 Let $a_1=0$, $a_2=2$, $(\rho_0,u_0)=(1,1)$, $(\rho_1,u_1)=(0.8,0.9)$ and
 $(\rho_2,u_2)=(0.9,0.7)$. Two SDWs will interact in a point
 $(X,T)=(12.2364,12.8427)$ with such data.
 \eex

 \begin{table}[htbp]
 \caption{$\;\;(\rho_0,u_0)=(1,1),(\rho_1,u_1)=(0.8,0.9),(\rho_2,u_2)=(0.9,0.7)$}\label{ds3}
 \vspace*{-0.3cm}
 $$
 \begin{array}{c|c|c|c|c|c|c|c|c|c|}
 \gamma & \kappa & \eps  & \rho_\eps &  u_\eps &   c_1   &   c_2   &  |Eq_1|         & |Eq_2| \\\hline
 2      & 0.5    & 0.5   & 1.16621   & 0.88674 & 0.20529 & 1.51813 &  8\cdot 10^{-5} & 4\cdot 10^{-5} \\\hline
 1.2    & 0.29   & 0.1   & 1.50419   & 0.86711 & 0.60356 & 1.11606 &  2\cdot 10^{-4} & 1\cdot 10^{-4} \\\hline
 1.02   & 0.099  & 0.01  & 3.75748   & 0.85659 & 0.80459 & 0.90592 &  4\cdot 10^{-3} & 2\cdot 10^{-3} \\\hline
 1.01   & 0.071  & 0.005 & 5.96447   & 0.85544 & 0.82632 & 0.88306 &  1\cdot 10^{-3} & 8\cdot 10^{-4} \\\hline
 1.006  & 0.055  & 0.003 & 8.66370   & 0.85341 & 0.83808 & 0.86341 &  2\cdot 10^{-2} & 3\cdot 10^{-2} \\\hline
 \end{array}
 $$
 \end{table}

 Neboj\v sa Dedovi\'c\\
 Department of Agricultural Engineering, University of Novi Sad\\
 Trg D.\ Obradovi\'{c}a 8, 21000 Novi Sad, Serbia\\
 dedovicn@uns.ac.rs\vspace*{0.1cm}\\

 Marko Nedeljkov\\
 Department of Mathematics and Informatics, University of Novi Sad\\
 Trg D.\ Obradovi\'{c}a 4, 21000 Novi Sad, Serbia\\
 markonne@uns.ac.rs

 \newpage

 \section{Appendix}

 %
 %

 %
 %

 %
 \begin{figure}[htbp]
 \centerline{
 \epsfig{figure=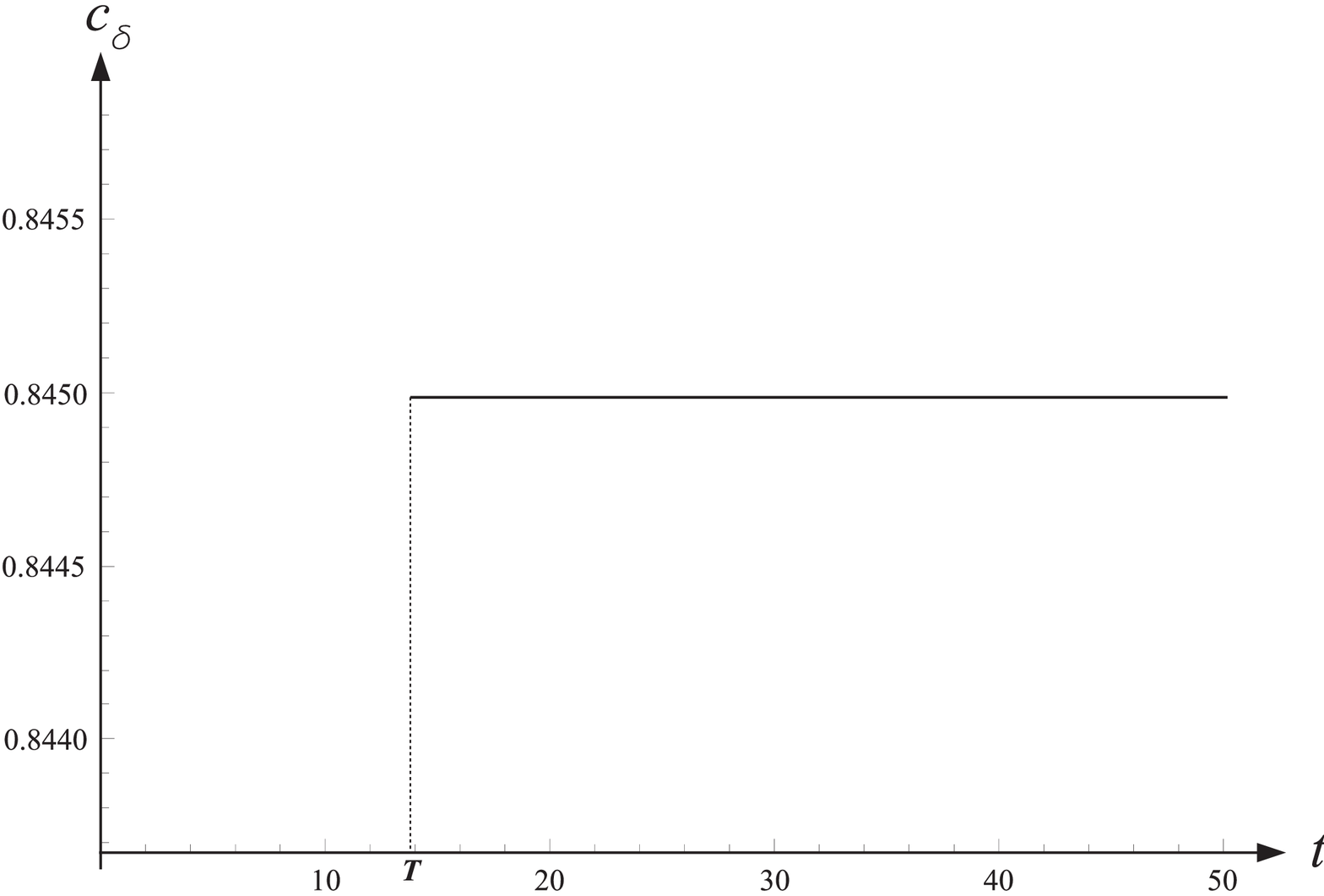,width=10cm
 }} \epsfxsize=1\textwidth
 \caption{Speed of delta shock formed after double delta shock
 interaction, Case A, Example \ref{ex1}.\label{usnew}}
 \end{figure}
 %

 %
 \begin{figure}[htbp]
 \centerline{
 \epsfig{figure=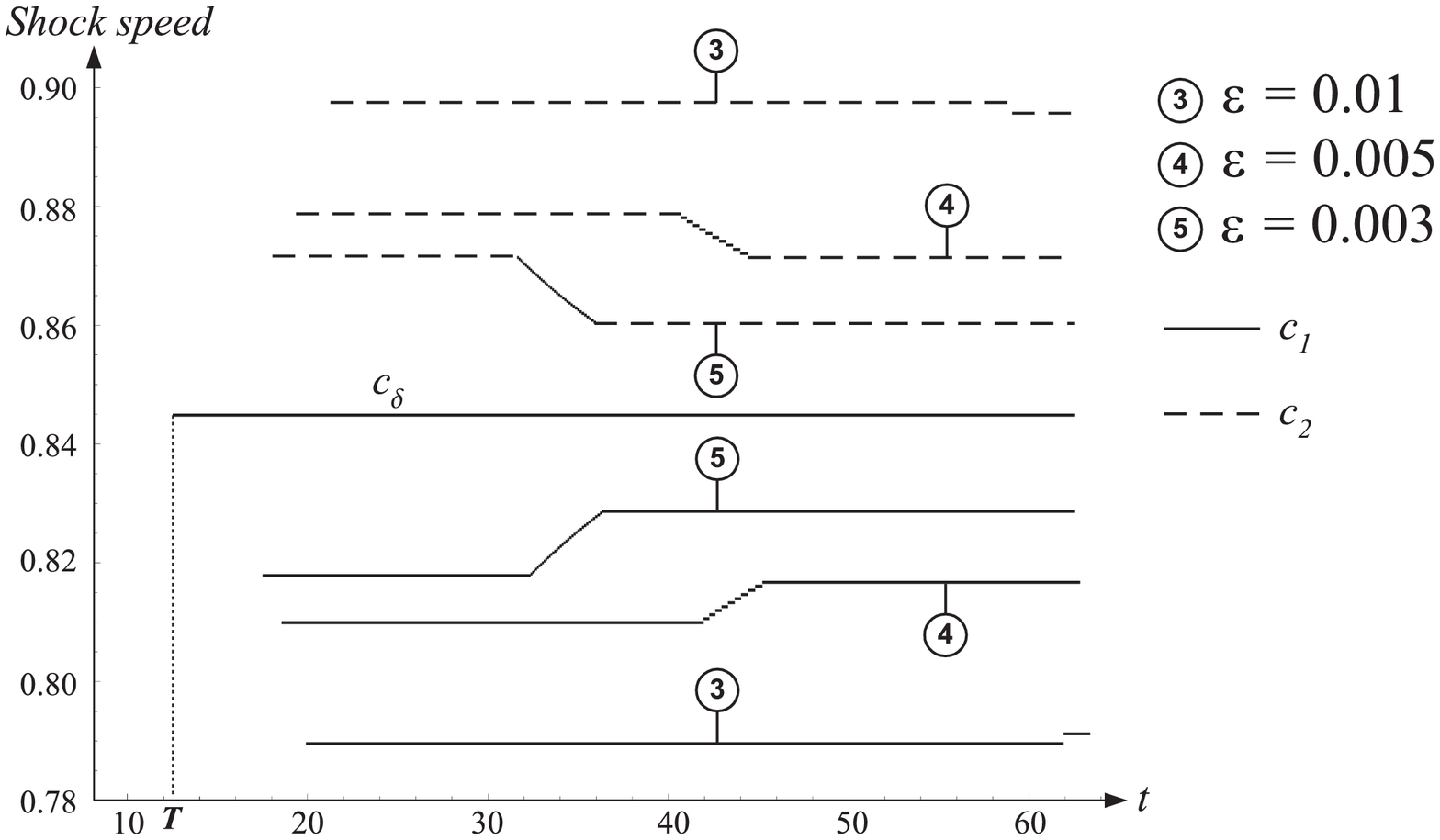,width=12cm
 }} \epsfxsize=1\textwidth
 \caption{Speed of the first ($S_1$) wave on the left hand side and the
 last ($S_2$) wave on the right hand side for various $\eps$, Case A, Example \ref{ex1}.\label{samobrzc}}
 \end{figure}
 %

 %
 \begin{figure}[htbp]
 \centerline{
 \epsfig{figure=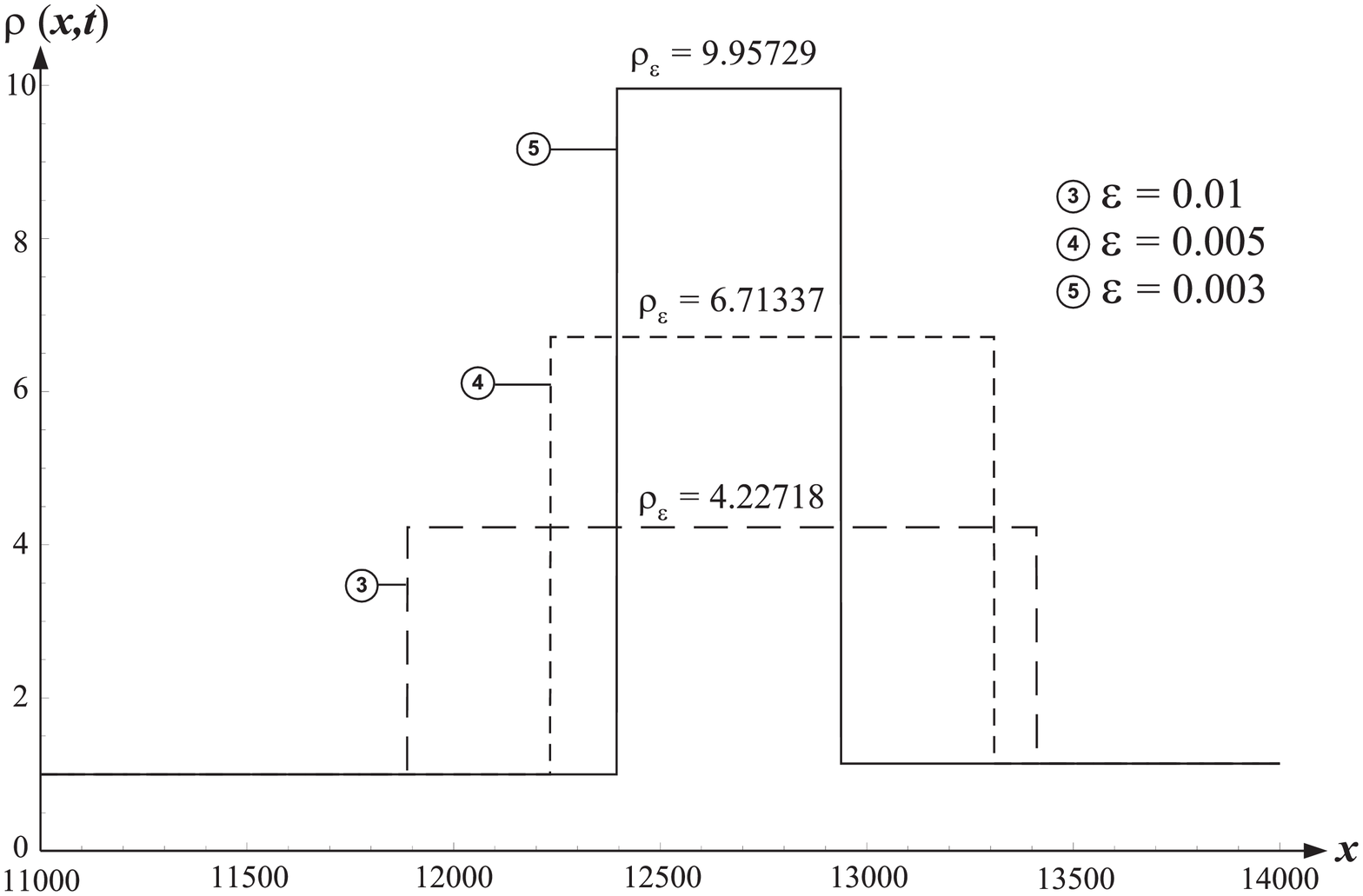,width=11cm
 }} \epsfxsize=1\textwidth
 \caption{Solution $\rho(x,t)$ for various $\eps$ at $t=15000$, Case A, Example \ref{ex1}.\label{us1newro}}
 \end{figure}
 %

 %
 \begin{figure}[htbp]
 \centerline{
 \epsfig{figure=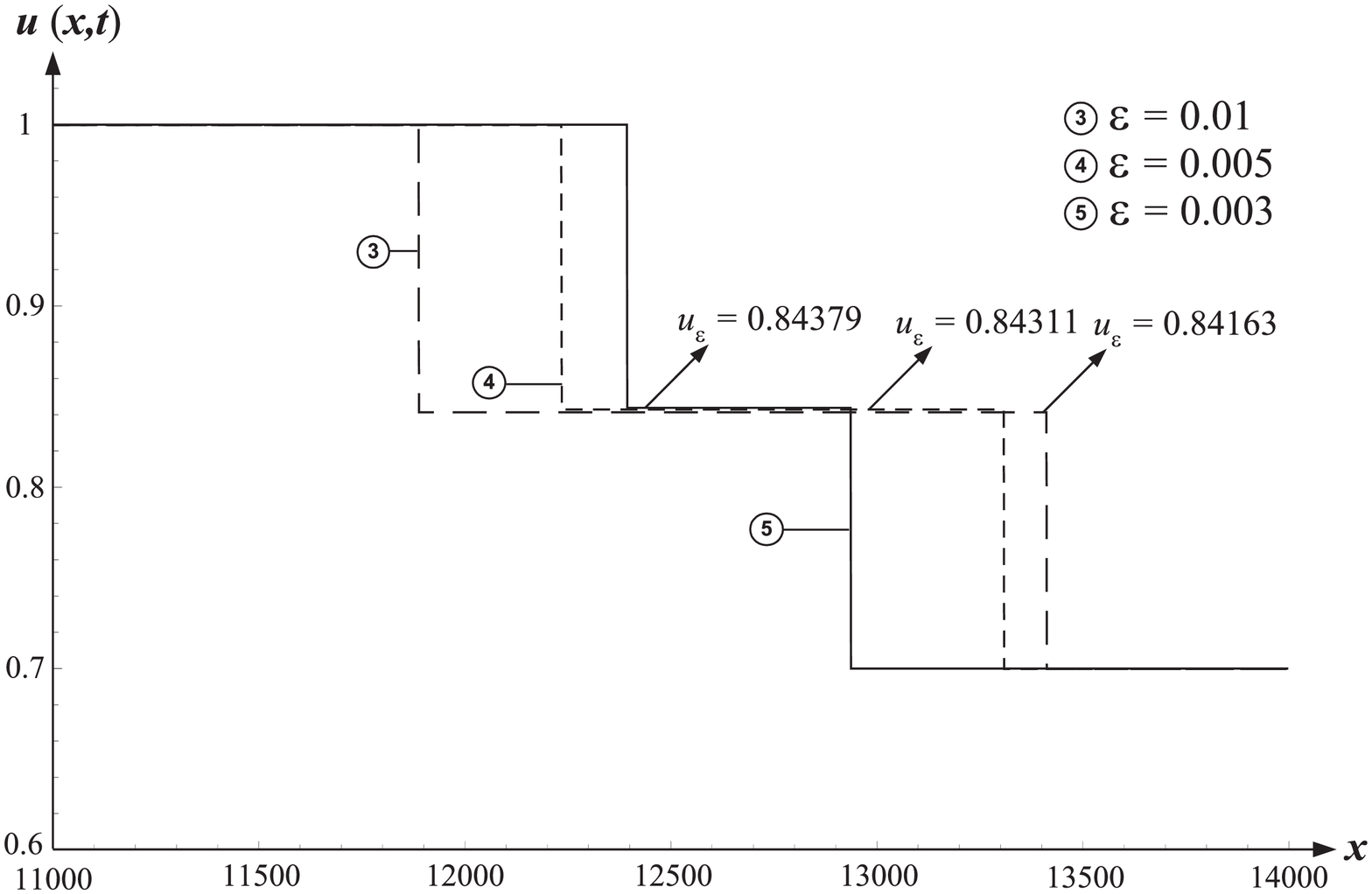,width=11cm
 }} \epsfxsize=1\textwidth
 \caption{Solution $u(x,t)$ for various $\eps$ at $t=15000$, Case A, Example \ref{ex1}.\label{us1newu}}
 \end{figure}
 %

 %
 %

 %
 \begin{figure}[htbp]
 \centerline{
 \epsfig{figure=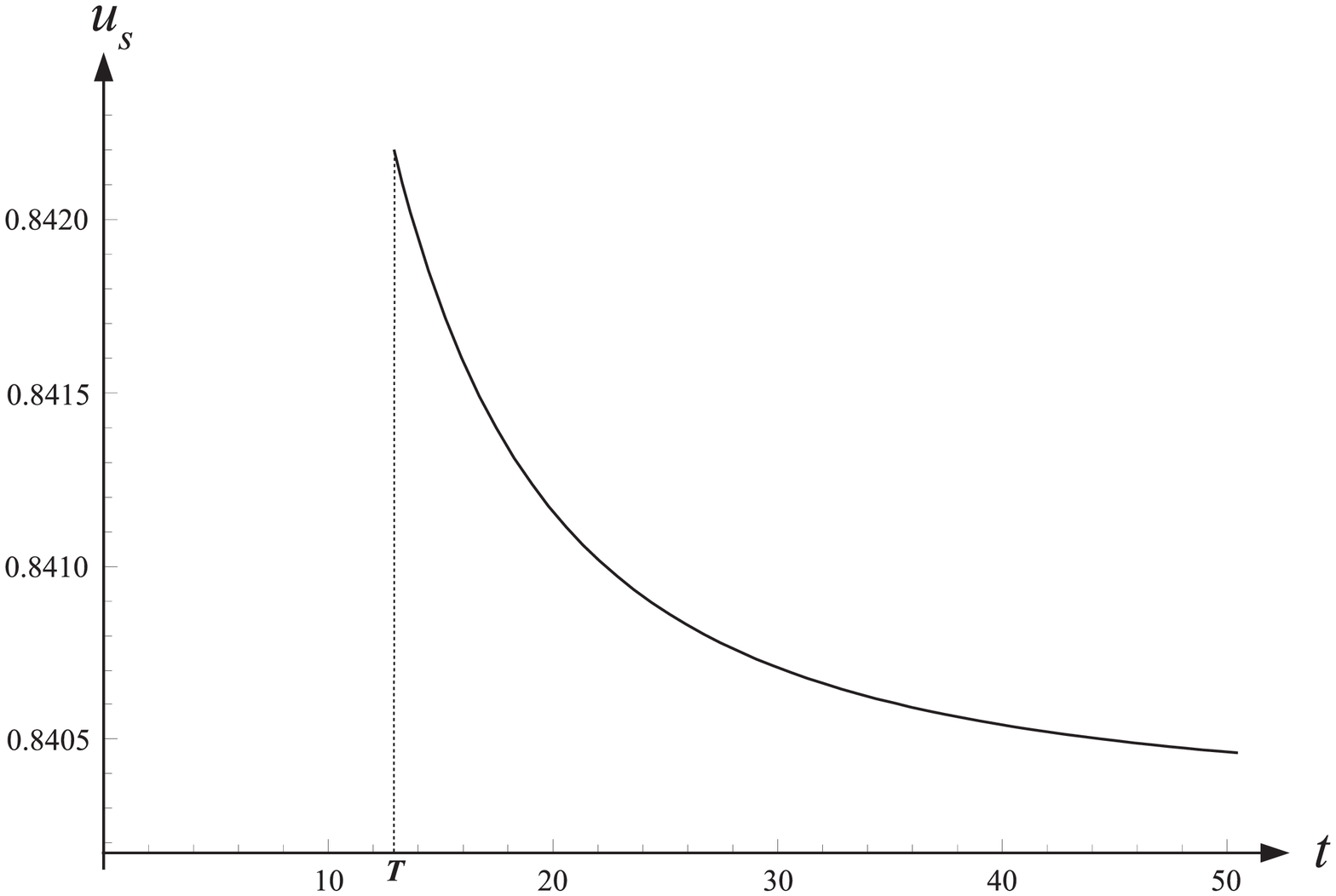,width=11cm
 }} \epsfxsize=1\textwidth
 \caption{Speed of delta shock formed after double delta shock
 interaction, Case B, Example \ref{ex2}.\label{us2newjul}}
 \end{figure}
 %

 %
 \begin{figure}[htbp]
 \centerline{
 \epsfig{figure=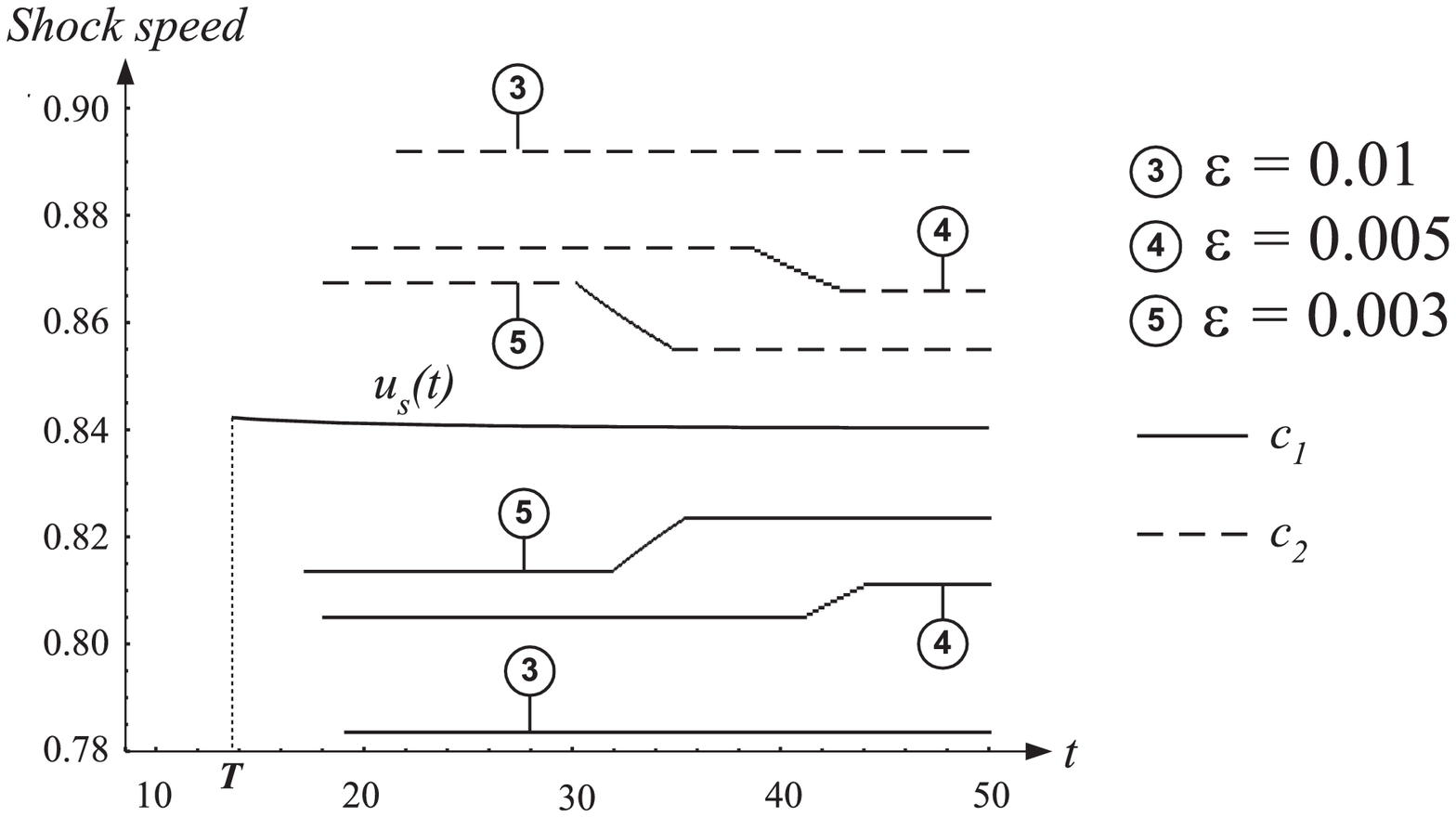,width=12cm
 }} \epsfxsize=1\textwidth
 \caption{Speed of the first ($S_1$) wave on the left hand side and the
 last ($S_2$) wave on the right hand side for various $\eps$, Case B, Example \ref{ex2}.\label{us2new}}
 \end{figure}
 %

 %
 \begin{figure}[htbp]
 \centerline{
 \epsfig{figure=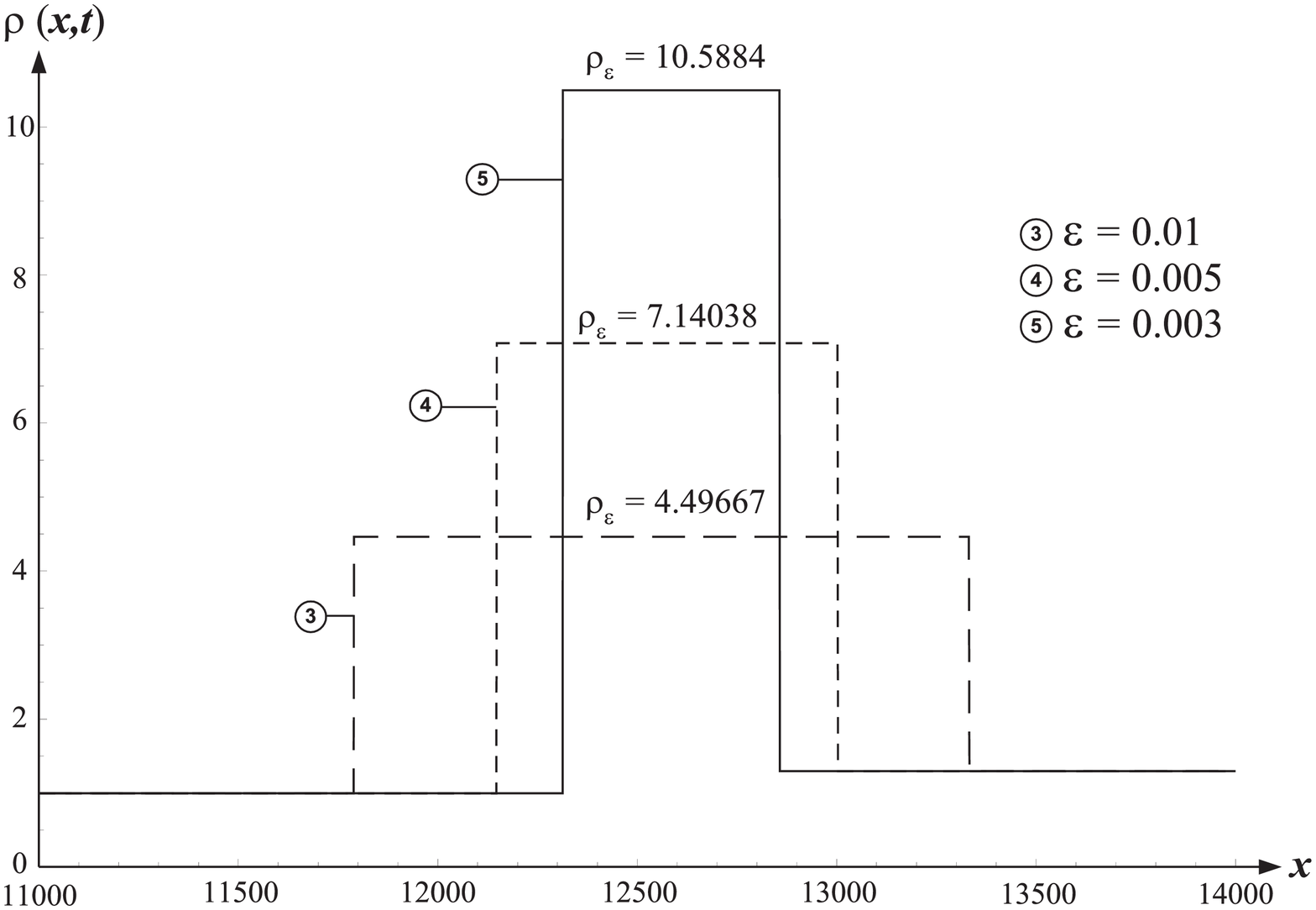,width=11cm
 }} \epsfxsize=1\textwidth
 \caption{Solution $\rho(x,t)$ for various $\eps$ at $t=15000$, Case B, Example \ref{ex2}.\label{us2newro}}
 \end{figure}
 %

 %
 \begin{figure}[htbp]
 \centerline{
 \epsfig{figure=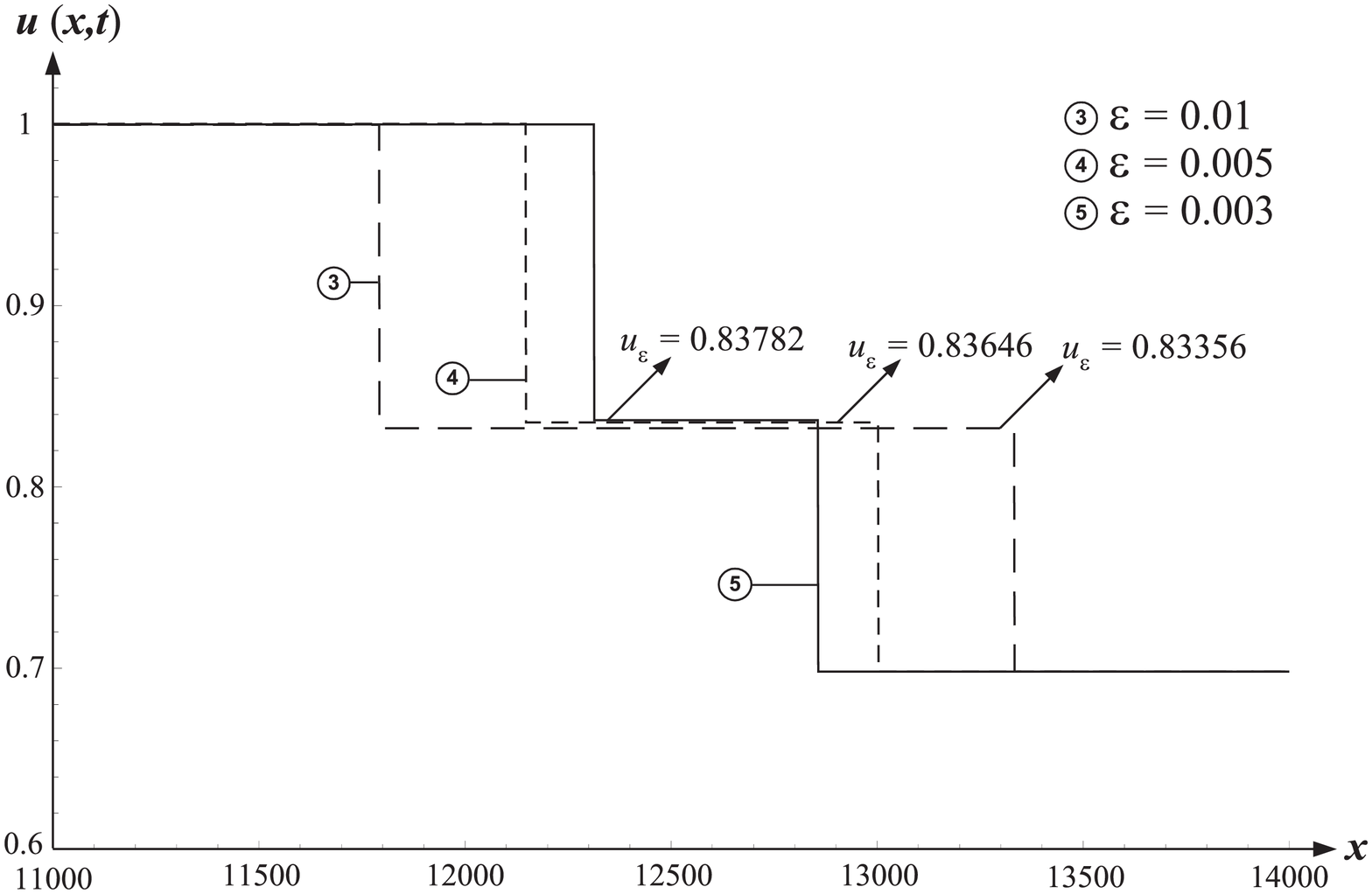,width=11cm
 }} \epsfxsize=1\textwidth
 \caption{Solution $u(x,t)$ for various $\eps$ at $t=15000$, Case B, Example \ref{ex2}.\label{us2newu}}
 \end{figure}
 %

 %
 %

 %
 \begin{figure}[htbp]
 \centerline{
 \epsfig{figure=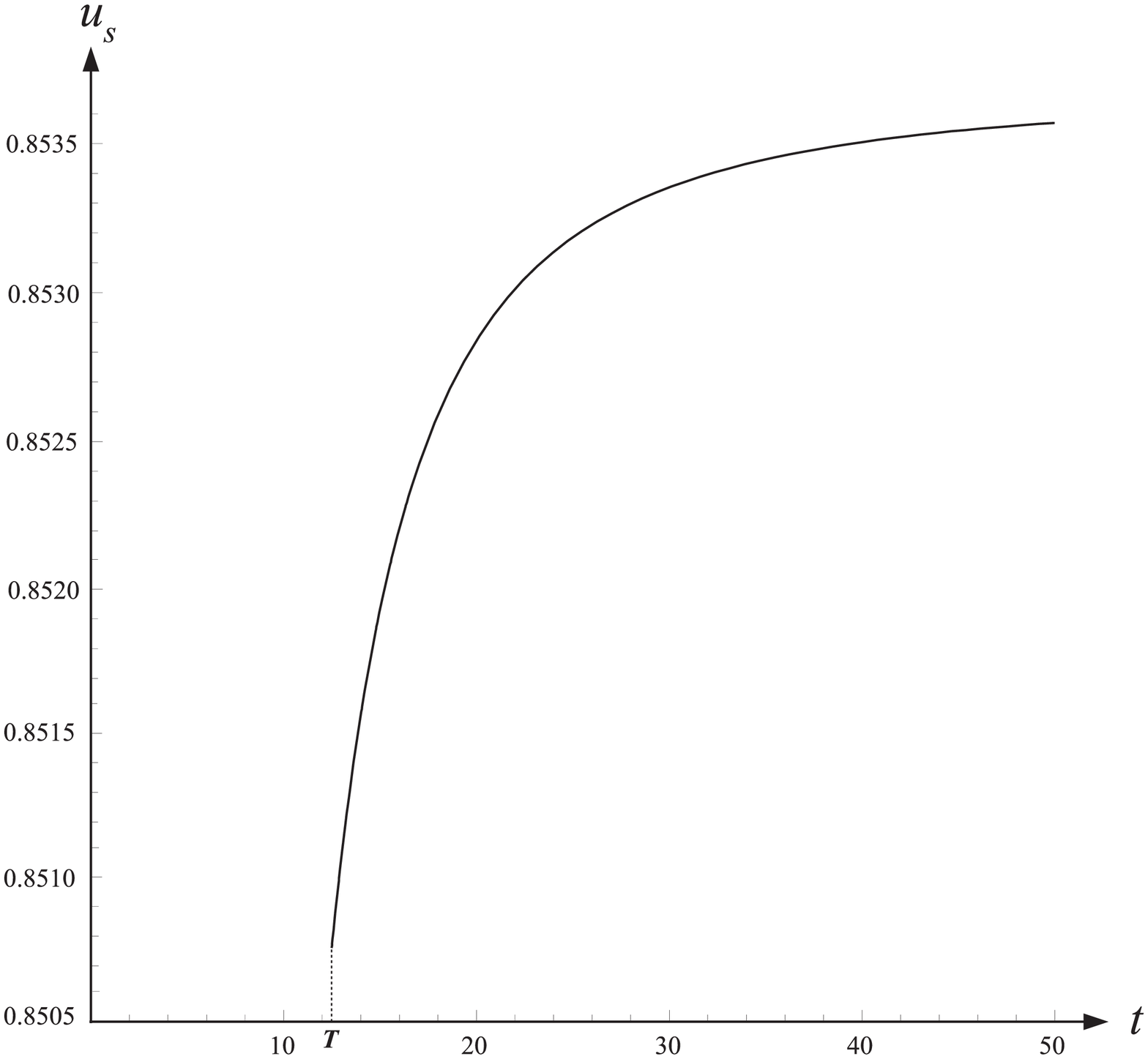,width=10cm
 }} \epsfxsize=1\textwidth
 \caption{Speed of delta shock formed after double delta shock interaction, Case B, Example \ref{ex3}.\label{us3newjul}}
 \end{figure}
 %

 %
 \begin{figure}[htbp]
 \centerline{
 \epsfig{figure=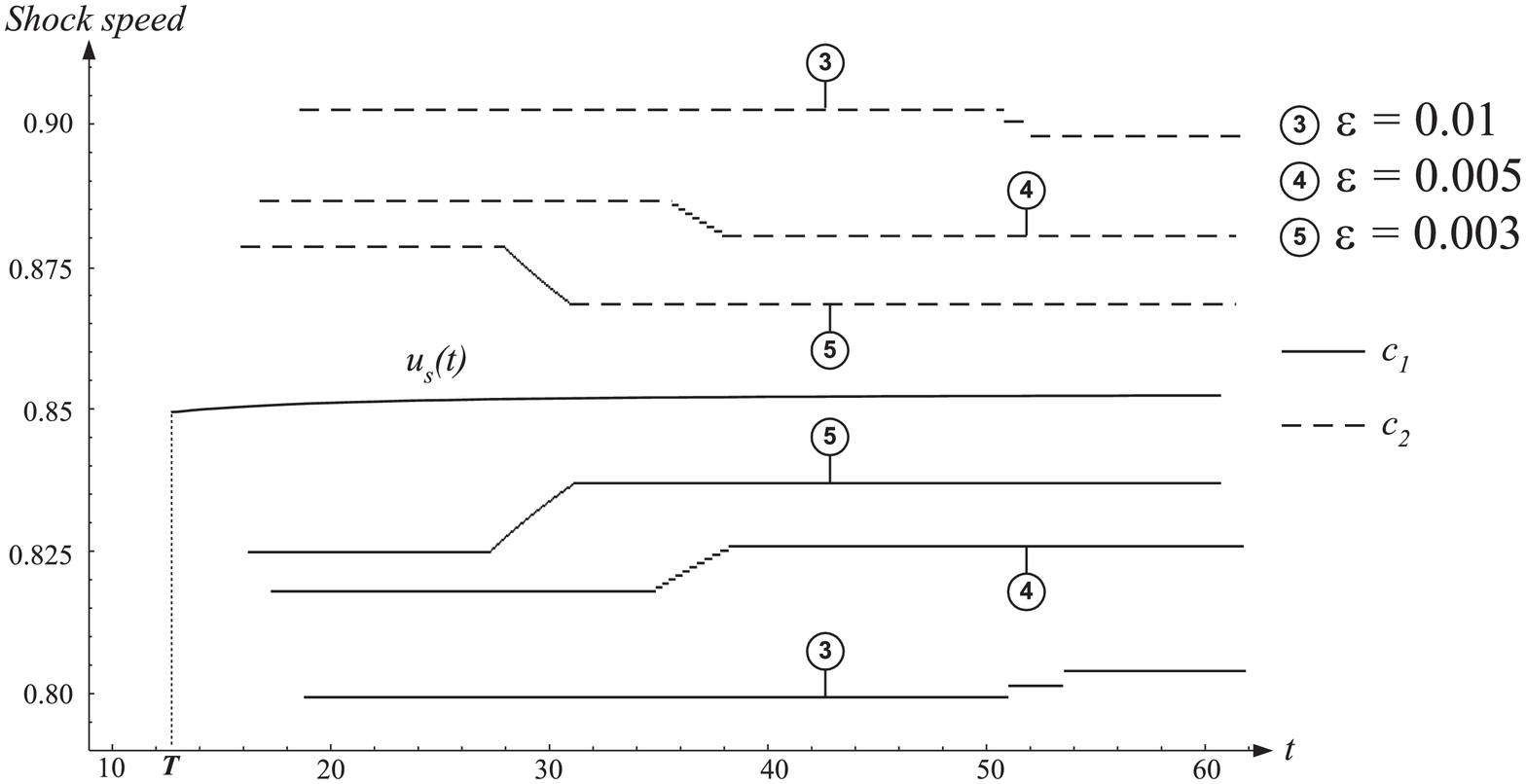,width=13cm
 }} \epsfxsize=1\textwidth
 \caption{Speed of the first ($S_1$) wave on the left hand side and the
 last ($S_2$) wave on the right hand side for various $\eps$, Case B, Example \ref{ex3}.\label{us3new}}
 \end{figure}
 %

 %
 \begin{figure}[htbp]
 \centerline{
 \epsfig{figure=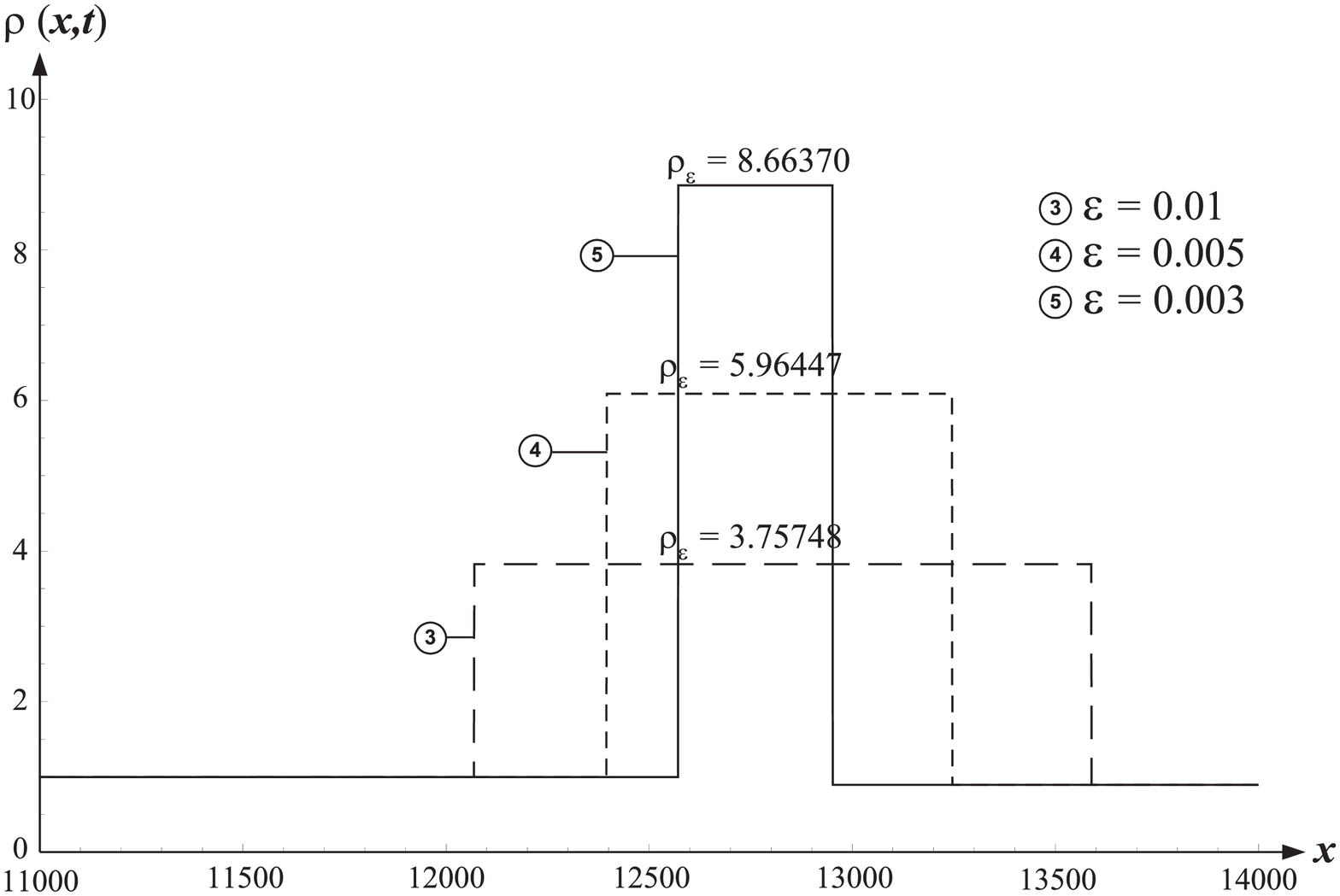,width=11cm
 }} \epsfxsize=1\textwidth
 \caption{Solution $\rho(x,t)$ for various $\eps$ at $t=15000$, Case B, Example \ref{ex3}.\label{us3newro}}
 \end{figure}
 %

 %
 \begin{figure}[htbp]
 \centerline{
 \epsfig{figure=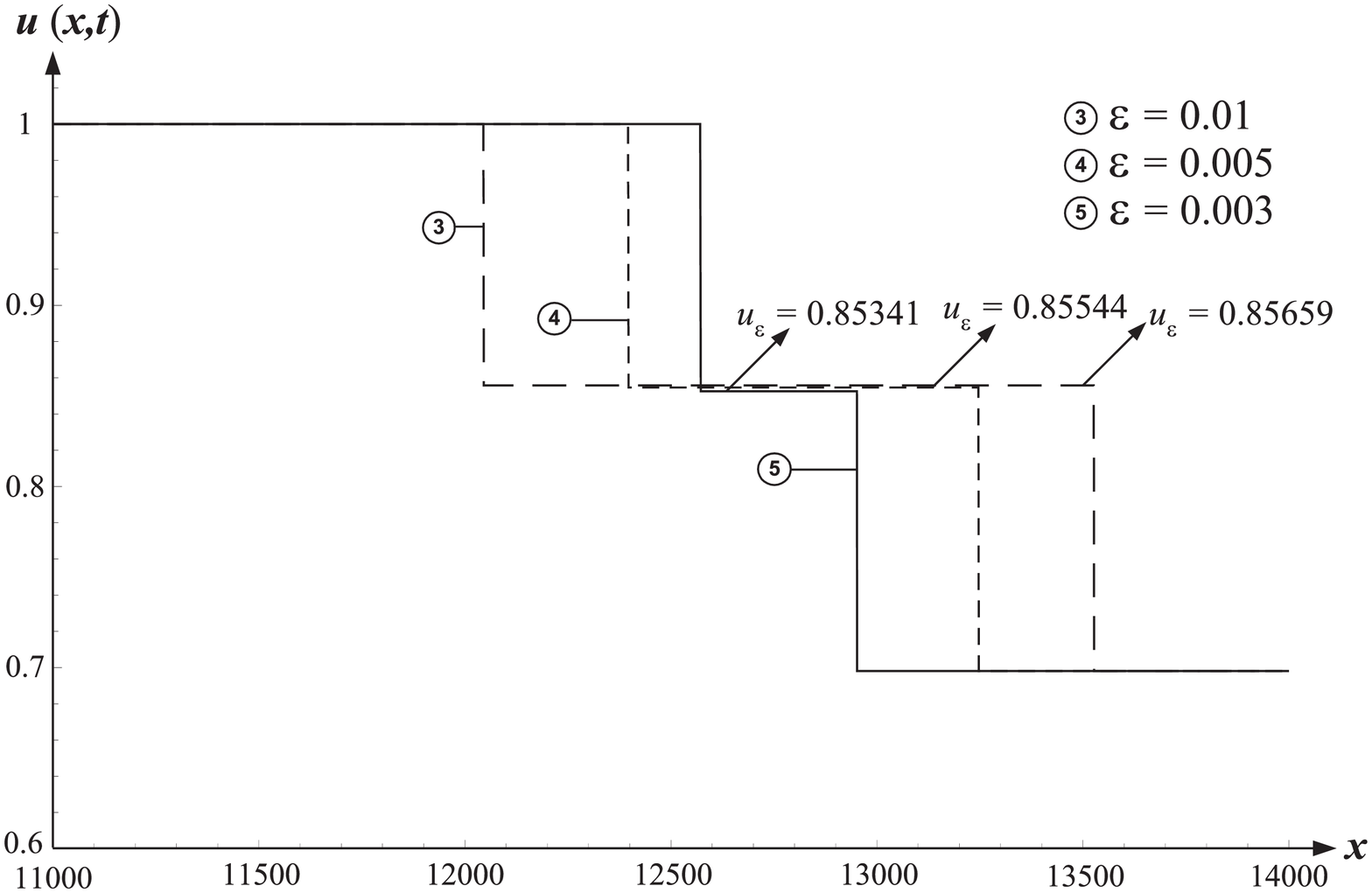,width=11cm
 }} \epsfxsize=1\textwidth
 \caption{Solution $u(x,t)$ for various $\eps$ at $t=15000$, Case B, Example \ref{ex3}.\label{us3newu}}
 \end{figure}
 %


\addcontentsline{toc}{part}{Bibliography}
 \begin{thebibliography}{9999}
 \bibitem{FA} F. Asakura, {\it Wave-front tracking method for the
 equations of isentropic gas dynamics},  Quart. Appl. Math. 63 (2005), no. 1, 20--33.
 %
 %
 %
 %
 %
 %
 \bibitem{AB} A. Bressan, {\it Hyperbolic Systems of Conservation Laws. The
 One-Dimensi-onal Cauchy Problem}, Oxford University Press, New York, 2000.
 %
 %
 %
 %
 %
 %
 %
 \bibitem{GQCHL} G.Q. Chen and H. Liu, {\it Formation of delta-shocks and vacuum states
 in the vanishing pressure limit of solutions to the isentropic Euler
 equations}, SIAM J. Math. Anal. 34 (2003), no. 4, 925--938.
 %
 %
 %
 %
 %
 %
 \bibitem{RCKF} R. Courant and K. Friedrichs,
 {\it Supersonic Flow and Shock Waves}, Applied Mathematical Sciences, Vol. 21. Springer-Verlag, New
York-Heidelberg, 1976.
 %
 %
 %
 %
 %
 %
 %
 %
 %
 %
 %
 %
 %
 %
 %
 %
 %
 %
 %
 %
 %
 %
 %
 %
 %
 %
 %
 %
 %
 %
 %
 %
 %
 %
 %
 %
 %
 %
 %
 %
 %
 %
 %
 %
 %
 %
 %
 %
 \bibitem{RL2} R.J. LeVeque, {\it Numerical Methods for Conservation Laws},
 Birkhäuser Verlag, Basel, 1990.
 %
 %
 %
 %
 %
 %
 %
 %
 %
 %
 \bibitem{DMMN} D. Mitrovi\'c and M. Nedeljkov, {\it Delta shock waves as a limit of shock waves},
  J. Hyperbolic Differ. Equ. 4 (2007), no. 4, 629--653.
 %
 %
 \bibitem{MN2} M. Nedeljkov, {\it Delta and singular delta locus for one
 dimensional systems of conservation laws},  Math. Methods Appl. Sci. 27 (2004), no. 8,
 931--955.
 %
 \bibitem{MN4} M. Nedeljkov, {\it Singular shock waves in interactions},
 Quart. Appl. Math. 66 (2008), no. 2, 281--302.
 %
 %
 \bibitem{MN2009p} M. Nedeljkov, Shadow waves -- entropies and interactions
 for delta and singular shocks (2009), to appear in Arch. Ration.
 Mech. Anal.
 %
 \bibitem{MN5} M. Nedeljkov, M. Oberguggengberger, {\it Interactions of delta shock waves in a
 strictly hyperbolic system of conservation laws}, J. Math. Anal. Appl. 344 (2008), no. 2, 1143--1157.
 %
 %
 \bibitem{TNJS} T. Nishida, J.A. Smoller,
 {\it Solutions in the Large for Some Nonlinear Hyperbolic Conservation Laws},
 Comm. Pure Appl. Math. 26 (1973), 183--200.
 %
 %
 %
 %
 %
 %
 %
 %
 %
 %
 %
 \bibitem{BR} B. Riemann, {\it Ueber die Fortpflanzung ebener Luftwellen
 von endlicher Schwingungsweite}, Gott.Abh.Math.Cl. 8 (1860), 43--65.
 %
 %
 %
 %
 \bibitem{ERS} E. Weinan, Y.G. Rykov, Ya.G. Sinai, {\it Generalized variotional principles,
 global weak solutions and behavior with random initial data for systems of consevation laws
 arising in adhesion particle dynamics},  Comm. Math. Phys. 177 (1996), no. 2, 349--380.
 %
 %
 %
 %
 %
 \bibitem{JS} J. Smoller, {\it Shock Waves and Reaction-Diffusion
 Equations}, Springer-Verlag, New York, 1994.
 %
 %
 %
 %
 %
 %
 %
 %
 %
 %
 \vspace*{0.6cm}
 \end{thebibliography}
 \end{document}